\font\eusm=eusm10                   


\font\eusms=eusm7                       

\font\eusmss=eusm5                      

\input amstex

\documentstyle{amsppt}
  \magnification=1100
  \hsize=6.2truein
  \vsize=9.0truein
  \hoffset 0.1truein
  \parindent=2em

\NoBlackBoxes

\newcount\mycitestyle \mycitestyle=1


\newcount\theTime
\newcount\theHour
\newcount\theMinute
\newcount\theMinuteTens
\newcount\theScratch
\theTime=\number\time
\theHour=\theTime
\divide\theHour by 60
\theScratch=\theHour
\multiply\theScratch by 60
\theMinute=\theTime
\advance\theMinute by -\theScratch
\theMinuteTens=\theMinute
\divide\theMinuteTens by 10
\theScratch=\theMinuteTens
\multiply\theScratch by 10
\advance\theMinute by -\theScratch
\def\timeHHMM{{\number\theHour:\number\theMinuteTens\number\theMinute}}

\def\today{{\number\day\space
 \ifcase\month\or
  January\or February\or March\or April\or May\or June\or
  July\or August\or September\or October\or November\or December\fi
 \space\number\year}}

\def\timeanddate{{\timeHHMM\space o'clock, \today}}



\define\Afr{{\frak A}}                       


\define\ah{{\hat a}}                         






\define\Ao{{A\oup}}                          








\define\biggnm#1{                            
  \bigg|\bigg|#1\bigg|\bigg|}

\define\bignm#1{                             
  \big|\big|#1\big|\big|}

\define\bh{\hat b}                           

\define\Bo{{B\oup}}                          






\define\clspan{\overline\lspan}              



\define\Cpx{\bold C}                         


\define\dif{\text{\it d}}                    

\define\dist{{\text{\rm dist}}}              






\define\eqdef{{\;\overset\text{def}\to=\;}}     







\define\fpamalg#1{{\dsize\;                  
     \operatornamewithlimits*_{#1}\;}}


\define\freeprod#1#2{\mathchoice             
     {\operatornamewithlimits{\ast}_{#1}^{#2}}
     {\raise.5ex\hbox{$\dsize\operatornamewithlimits{\ast}
      _{#1}^{#2}$}\,}
     {\text{oops!}}{\text{oops!}}}

\define\freeprodi{\mathchoice                
     {\operatornamewithlimits{\ast}
      _{\iota\in I}}
     {\raise.5ex\hbox{$\dsize\operatornamewithlimits{\ast}
      _{\sssize\iota\in I}$}\,}
     {\text{oops!}}{\text{oops!}}}









\define\GL{{\text{\rm GL}}}                  












\define\Ic{{\Cal I}}                         








\define\Integers{\bold Z}                    


\define\KHil{{\mathchoice                    
     {\text{\eusm K}}
     {\text{\eusm K}}
     {\text{\eusms K}}
     {\text{\eusmss K}}}}







\define\lrnm#1{\left|\left|#1\right|\right|} 

\define\lspan{\text{\rm span}@,@,@,}         

\define\Mfr{{\frak M}}                       








\define\nm#1{||#1||}                         


\define\Naturals{{\bold N}}                  




\define\onehat{\hat 1}                       



\define\oup{^{\text{\rm o}}}                 

\define\owedge{{                             
     \operatorname{\raise.5ex\hbox{\text{$
     \ssize{\,\bigcirc\llap{$\ssize\wedge\,$}\,}$}}}}}

\define\owedgeo#1{{                          
     \underset{\raise.5ex\hbox
     {\text{$\ssize#1$}}}\to\owedge}}










\define\QED{\newline                         
            \line{$\hfill$\qed}\enddemo}









\define\Reals{{\bold R}}                     


\define\restrict{\lower .3ex                 
     \hbox{\text{$|$}}}




\define\smd#1#2{\underset{#2}\to{#1}}          

\define\smdb#1#2{\undersetbrace{#2}\to{#1}}    

\define\smdbp#1#2#3{\overset{#3}\to            
     {\smd{#1}{#2}}}

\define\smdbpb#1#2#3{\oversetbrace{#3}\to      
     {\smdb{#1}{#2}}}

\define\smdp#1#2#3{\overset{#3}\to             
     {\smd{#1}{#2}}}

\define\smdpb#1#2#3{\oversetbrace{#3}\to       
     {\smd{#1}{#2}}}

\define\smp#1#2{\overset{#2}\to                
     {#1}}                                     

\define\solim{{\text{\rm s.o.--}\lim}}       



\define\sr{{\text{\rm sr}}}                  




\define\Tcirc{\bold T}                       

\define\tint{{\tsize\int}}                   

\define\tocdots                              
  {\leaders\hbox to 1em{\hss.\hss}\hfill}    


\define\Tr{\text{\rm Tr}}                    


\define\Uc{{\Cal U}}                         


\define\uh{{\hat u}}                         























  \newcount\bibno \bibno=0
  \def\newbib#1{\advance\bibno by 1 \edef#1{\number\bibno}}
  \ifnum\mycitestyle=1 \def\cite#1{{\rm[\bf #1\rm]}} \fi
  \def\scite#1#2{{\rm[\bf #1\rm, #2]}}


  \newcount\ignorsec \ignorsec=0
  \def\notasec{\ignorsec=1}
  \def\forgetsec{\ignorsec=1}

  \newcount\secno \secno=0
  \def\newsec#1{\procno=0 \subsecno=0 \ignorsec=0
    \advance\secno by 1 \edef#1{\number\secno}
    \edef\currentsec{\number\secno}}

  \newcount\subsecno
  \def\newsubsec#1{\procno=0 \advance\subsecno by 1 \edef#1{\number\subsecno}
    \edef\currentsec{\number\secno.\number\subsecno}}

  \newcount\appendixno \appendixno=0
  \def\newappendix#1{\procno=0 \ignorsec=0 \advance\appendixno by 1
    \ifnum\appendixno=1 \edef\appendixalpha{\hbox{A}}
      \else \ifnum\appendixno=2 \edef\appendixalpha{\hbox{B}} \fi
      \else \ifnum\appendixno=3 \edef\appendixalpha{\hbox{C}} \fi
      \else \ifnum\appendixno=4 \edef\appendixalpha{\hbox{D}} \fi
      \else \ifnum\appendixno=5 \edef\appendixalpha{\hbox{E}} \fi
      \else \ifnum\appendixno=6 \edef\appendixalpha{\hbox{F}} \fi
    \fi
    \edef#1{\appendixalpha}
    \edef\currentsec{\appendixalpha}}

  \newcount\procno \procno=0
  \def\newproc#1{\advance\procno by 1
   \ifnum\ignorsec=0 \edef#1{\currentsec.\number\procno}
   \else \edef#1{\number\procno}
   \fi}

  \newcount\tagno \tagno=0
  \def\newtag#1{\advance\tagno by 1 \edef#1{\number\tagno}}



\notasec
  \newtag{\RedFreeProd}
  \newtag{\AvitzourCond}
  \newtag{\AfrfdAbelian}
  \newtag{\SimplicityCond}
\newsec{\StmntMainResults}
\forgetsec
  \newtag{\ABfdAbelian}
 \newproc{\FPfdAbelian}
  \newtag{\AfrMainRes}
 \newproc{\FPdas}
\newsec{\Preliminaries}
 \newproc{\DefDiffuse}
 \newproc{\DASHU}
 \newproc{\TensorMn}
 \newproc{\ShortExactSeq}
 \newproc{\AheredBsrone}
 \newproc{\FullSimple}
 \newproc{\TwoProj}
 \newproc{\FreeProdDirectSum}
 \newproc{\notCpx}
\newsec{\InPresOfDiff}
 \newproc{\TechnicalLemma}
  \newtag{\xoneab}
  \newtag{\uNBAB}
  \newtag{\limLtwo}
  \newtag{\zcongxzero}
  \newtag{\usualterms}
  \newtag{\urxus}
 \newproc{\FPHaarUnitaryDixmier}
  \newtag{\dixmierprop}
  \newtag{\DixmierFive}
 \newproc{\TechnicalLemmatwo}
  \newtag{\zxz}
 \newproc{\FPHaarUnitarySRone}
 \newproc{\ExistsDAS}
 \newproc{\ExistsDAScorr}
\newsec{\FPFinManyS}
 \newproc{\AzeroCCC}
  \newtag{\NotationAzeroCdC}
  \newtag{\AzeroCCCcases}
  \newtag{\CaseIIAone}
  \newtag{\CaseIIdenseker}
  \newtag{\CaseIVAone}
  \newtag{\CaseIVdenseker}
  \newtag{\CaseVIAone}
  \newtag{\CaseVIdenseker}
 \newproc{\AzeroCCCC}
  \newtag{\LplusLzero}
  \newtag{\AzeroCCCCAfrone}
  \newtag{\AzeroCCCCAfr}
  \newtag{\CaseIIAzeroCCCCdense}
  \newtag{\CaseIIAzeroCCCCfullin}
  \newtag{\CaseIIAzeroCCCCfullintwo}
  \newtag{\CaseIIIAzeroCCCCfullini}
  \newtag{\CaseIIIAzeroCCCCfullinj}
 \newproc{\CdCCC}
  \newtag{\LplusLzeron}
  \newtag{\CdCCCdensei}
  \newtag{\CdCCCdensej}
  \newtag{\Lvarious}
  \newtag{\CdCCCrefpt}
  \newtag{\CaseIICdCCCdense}
  \newtag{\Aonezerocapj}
  \newtag{\Atwozerocapj}
 \newproc{\AzeroCdCCC}
  \newtag{\AzeroCdCCCAfr}
  \newtag{\AzeroCdCCCdensei}
  \newtag{\AzeroCdCCCdensej}
  \newtag{\AzeroCdCCCdense}
  \newtag{\AzeroCdCCCAfronefull}
 \newproc{\CdCCdC}
  \newtag{\CdCCdCdense}
 \newproc{\NonUnitalToo}
 \newproc{\AzeroCdCCdC}
  \newtag{\eqAzeroCdCCdC}
  \newtag{\eqAzeroCdCBzeroCdC}
 \newproc{\FinManyFDAbel}
  \newtag{\Aiota}
  \newtag{\LpluszeroMTTFDA}
 \newproc{\FinManyAzeroCdC}
  \newtag{\AzeroCdC}
 \newproc{\FinManyCofXAtoms}
  \newtag{\AfrCofXIsol}
\newsec{\CofXS}
 \newproc{\InvLimIsolAtomsDef}
 \newproc{\InvLimIsolAtomsEx}
 \newproc{\InvLimIsolAtomsFP}
  \newtag{\AfrIndLim}
  \newtag{\lambdanatoms}
\newsec{\FPInfManyS}
 \newproc{\FPInfMany}
\newsec{\Conjectures}
 \newproc{\ConjCofX}
 \newproc{\CofXNec}
 \newproc{\ConjFD}
  \newtag{\ConjFDFP}
  \newtag{\ConjFDnj}
  \newtag{\ConjFDmj}

\newbib{\AndersonBlackadarHaagerup}
\newbib{\Avitzour}
\newbib{\BrownLGZZStIsoHered}
\newbib{\BrownGreenRieffel}
\newbib{\DykemaZZFreeDim}
\newbib{\DykemaZZFPFDNT}
\newbib{\DykemaZZFaithful}
\newbib{\DykemaZZSNU}
\newbib{\DykemaHaagerupRordam}
\newbib{\DykemaRordamZZPI}
\newbib{\GermainZZKKeq}
\newbib{\GermainZZKthFFP}
\newbib{\HermanVassersteinZZsr}
\newbib{\PaschkeSalinasZZFreeProdGp}
\newbib{\PowersZZCsFtwo}
\newbib{\RieffelZZMoritaEqOpAlg}
\newbib{\RieffelZZsr}
\newbib{\RordamZZUnitaryRank}
\newbib{\VoiculescuZZSymmetries}
\newbib{\VoiculescuZZMult}
\newbib{\VDNbook}

\topmatter

  \title Simplicity and the stable rank of some free product
         C$^*$--algebras
  \endtitle

  \leftheadtext{Simplicity of free products \timeanddate}
  \rightheadtext{Simplicity of free products \timeanddate}

  \author Kenneth J\. Dykema
  \endauthor

  \date \timeanddate \enddate

  \affil 
    Department of Mathematics and Computer Science \\
    Odense Universitet, Campusvej 55 \\
    DK-5230 Odense M \\
    Denmark
  \endaffil

  \address
    Department of Mathematics and Computer Science,
    Odense Universitet, Campusvej 55,
    DK-5230 Odense M,
    Denmark
  \endaddress

  \abstract
   A necessary and sufficient condition for the simplicity of the
   C$^*$--algebra reduced free product
   of finite dimensional abelian algebras is found, and it is proved that the
   stable rank of every such free product is~1.
   Related results about other reduced free products of C$^*$--algebras
   are proved.
  \endabstract


\endtopmatter

\document \TagsOnRight \baselineskip=18pt

\vskip3ex
\noindent{\bf Introduction.}
\vskip3ex

The reduced free product of C$^*$--algebras with respect to given states was
introduced independently by Voiculescu~\cite{\VoiculescuZZSymmetries} and
Avitzour~\cite{\Avitzour}.
It is the appropriate construction associated to Voiculescu's free
probability theory, (see~\cite{\VoiculescuZZSymmetries},~\cite{\VDNbook}).
The motivating example
concerns reduced group
C$^*$--algebras.
For a discrete group $G$, its reduced group C$^*$--algebra is
generated by the left regular representation of $G$ on $l^2(G)$
and is denoted $C^*_r(G)$.
Its canonical tracial state (which is the vector state associated to the
characteristic function of the identity element of $G$) is written $\tau_G$.
Then for discrete groups $G_1$ and $G_2$, the reduced free product construction
yields
$$ (C^*_r(G_1),\tau_{G_1})*(C^*_r(G_2),\tau_{G_2})=(C^*_r(G),\tau_G), $$
where $G=G_1*G_2$ is the free product of groups.

Voiculescu's definition of freeness is an abstraction of some essential
facets of the relationship between the copies
of $C^*_r(G_1)$ and $C^*_r(G_2)$ embedded in $C^*_r(G)$, with respect to the
trace $\tau_G$.
The reduced free product of C$^*$--algebras can be described with respect to
freeness as follows.
Let $A_1$ and $A_2$ be unital C$^*$--algebras with states $\phi_1$ and
respectively $\phi_2$, whose associated GNS representations are faithful.
Then the reduced free product of $(A_1,\phi_1)$ and $(A_2,\phi_2)$ is the
(unique) unital C$^*$--algebra $\Afr$ and state $\phi$ with unital embeddings
$A_\iota\hookrightarrow\Afr$ such that
\roster
\item the GNS representation associated to $\phi$ is faithful on $\Afr$;
\item $\phi\restrict_{A_\iota}=\phi_\iota$;
\item $A_1$ and $A_2$ are free with respect to $\phi$;
\item $\Afr$ is generated by $A_1\cup A_2$.
\endroster
We denote this by
$$ (\Afr,\phi)=(A_1,\phi_1)*(A_2,\phi_2). \tag{\RedFreeProd} $$
It is further known~\cite{\VoiculescuZZSymmetries} (or
see~\scite{\VDNbook}{2.5.3})
that $\phi$ is a trace if $\phi_1$ and $\phi_2$ are traces.
Moreover, by~\cite{\DykemaZZFaithful}, $\phi$ is also faithful on $\Afr$ if
$\phi_1$ and $\phi_2$ are faithful.

The reduced free product thus provides a multitude of constructions of
C$^*$--algebras, about which some results are known
(see~\cite{\VoiculescuZZSymmetries}, \cite{\Avitzour},
\cite{\GermainZZKKeq},
\cite{\DykemaRordamZZPI}, \cite{\DykemaHaagerupRordam}).
For example, many can be distinguished one from the other using K--theory,
(see~\cite{\GermainZZKKeq} and~\cite{\GermainZZKthFFP}).
However, questions abound.

Perhaps the most basic question concerns simplicity of reduced free product
C$^*$--algebras.
In~\cite{\PowersZZCsFtwo}, R.T\. Powers showed that the reduced group
C$^*$--algebra
of the free group on two generators, $C^*_r(F_2)$, is simple and has unique
tracial state.
Paschke and Salinas~\cite{\PaschkeSalinasZZFreeProdGp} then proved the same
for
$C^*_r(G)$ whenever $G=G_1*G_2$ is the free product of groups, where $G_1$ has
at least two elements and $G_2$ has at least three.
Avitzour~\cite{\Avitzour} generalized further and showed that, for
the reduced free product~(\RedFreeProd), $\Afr$ is simple if there are
unitaries $u,v\in A_1$ and $w\in A_2$ such that
$$ \alignedat 2
\phi_1(u)=\phi_1(v)&=0=\phi_1(u^*v),\qquad&\phi_1(u^*\cdot u)&=\phi_1, \\
\phi_2(w)&=0,&\phi_2(w^*\cdot w)&=\phi_2.
\endalignedat \tag{\AvitzourCond} $$
(Actually, Avitzour required also $\phi_1$ and $\phi_2$ to be faithful, but
this hypothesis is easily dispensed with.)

Avitzour's conditions imply simplicity of many reduced free product
C$^*$--algebras, but there are plenty of cases where Avitzour's conditions are
not satisfied (see~\S4 of~\cite{\DykemaHaagerupRordam}), yet
intuition (or the analogous result for von Neumann algebras,
see~\cite{\DykemaZZFreeDim}) suggests the algebra is simple.
In this paper we give necessary and sufficient conditions for
simplicity of the reduced free product of arbitrary finite
dimensional abelian C$^*$--algebras.
Stated briefly, if
$$ (\Afr,\tau)=(A,\tau_A)*(B,\tau_B), \tag{\AfrfdAbelian} $$
where $A$ and $B$ are
finite dimensional abelian C$^*$--algebras satisfying $\dim(A)\ge3$ and
$\dim(B)\ge2$
and with faithful tracial states,
$\tau_A$ and $\tau_B$, then $\Afr$ is simple if and only if whenever $p$ is a
minimal projection of $A$ and $q$ is a minimal projection of $B$, we have
$$ \tau_A(p)+\tau_B(q)<1. \tag{\SimplicityCond} $$
The necessity of this condition can be seen
from~\cite{\AndersonBlackadarHaagerup}.
Note that the condition from~\cite{\DykemaZZFreeDim} for the analogous free
product of von Neumann algebras to be a factor is~(\SimplicityCond) but with
the strict inequality replaced by $\le$.
In addition, when the free product algebra, $\Afr$, from~(\AfrfdAbelian) is not
simple, our analysis allows one to easily find all ideals of $\Afr$.

We also show that for every reduced free product C$^*$--algebra, $\Afr$, as
in~(\AfrfdAbelian), the stable rank of~$\Afr$ is~1, regardless of the
simplicity of $\Afr$.
The (topological) stable rank was invented by M.A\.~Rieffel
in~\cite{\RieffelZZsr} as a sort of dimension for C$^*$--algebras, and was
in~\cite{\HermanVassersteinZZsr} shown to be equal to the Bass stable rank.
The first result about the stable rank of reduced free product C$^*$--algebras
is in~\cite{\DykemaHaagerupRordam}, where it is proved that the free product
with respect to traces of C$^*$--algebras $A_1$ and $A_2$ has stable rank~1 if
the Avitzour conditions~(\AvitzourCond) are satisfied.
Hence the present paper's results regarding stable rank are, for a restricted
class of C$^*$--algebra reduced free products, a considerable generalization
and lend support to the plausible conjecture that every reduced free product of
C$^*$--algebras with respect to faithful, tracial states has stable rank~1.

In~\S\StmntMainResults{} we state the main results proved in the paper.
In~\S\Preliminaries{} concepts and results essential for the sequel are
covered, including some dealing with stable rank, full hereditary subalgebras
and free products.
In~\S\InPresOfDiff{} we prove simplicity and stable rank~1
for free products of two C$^*$--algebras, when one is diffuse in a specific
sense.
In~\S\FPFinManyS{} the results about the free product of finite dimensional
abelian algebras are proved, as well as some related results about free
products of more general algebras.
In~\S\CofXS{} free products of abelian C$^*$--algebras with states that are
inductive limits of the algebras considered in~\S\FPFinManyS{} are proved.
In~\S\FPInfManyS{} we consider free products of infinitely many finite
dimensional abelian C$^*$--algebras.
Finally, in~\S\Conjectures{} we make two conjectures about simplicity of other
free product C$^*$--algebras.

\vskip3ex
\noindent{\bf\S\StmntMainResults. Statement of the main results.}
\nopagebreak
\vskip3ex
In this section, we state the main results in more detail.
Let
$$ \aligned
(A,\tau_A)&=\smdp\Cpx{\alpha_1}{p_1}\oplus\smdp\Cpx{\alpha_2}{p_2}
\oplus\cdots\oplus\smdp\Cpx{\alpha_n}{p_n} \\ \vspace{1ex}
(B,\tau_B)&=\smdp\Cpx{\beta_1}{q_1}\oplus\smdp\Cpx{\beta_2}{q_2}
\oplus\cdots\oplus\smdp\Cpx{\beta_m}{q_m}.
\endaligned \tag{\ABfdAbelian} $$
This notation means that $n\in\Naturals$, that
$$ A=\undersetbrace{n\text{ times}}\to{\Cpx\oplus\cdots\oplus\Cpx}, $$
that $p_k$ is the projection
$$ p_k=\undersetbrace{k-1}\to{0\oplus\cdots\oplus0}\oplus1\oplus
\undersetbrace{n-k}\to{0\oplus\cdots\oplus0}, $$
and that $\tau_A$ is the state on $A$ given by $\tau_A(p_k)=\alpha_k$.
We thus need $\alpha_k\ge0$ and $\sum_1^n\alpha_k=1$, and because we want the
GNS representation of $\tau_A$ to be faithful we will always take
$\alpha_k>0$.
The same considerations apply to $(B,\tau_B)$ in~(\ABfdAbelian).
\proclaim{Theorem \FPfdAbelian}
Let $(A,\tau_A)$ and $(B,\tau_B)$ be finite dimensional, abelian
C$^*$--algebras with faithful states as in~(\ABfdAbelian), with $\dim(A)\ge3$
and $\dim(B)\ge2$.
Let
$$ (\Afr,\tau)=(A,\tau_A)*(B,\tau_B) $$
be the reduced free product of C$^*$--algebras.
Let
$$ \align
L_+&=\{(i,j)\mid\alpha_i+\beta_j>1\}\\
L_0&=\{(i,j)\mid\alpha_i+\beta_j=1\}.
\endalign $$
Then the stable rank of $\Afr$ is~1 and
$$ \Afr=\smp{\Afr_0}{r_0}\oplus\bigoplus_{(i,j)\in L_+}
\smdp\Cpx{\alpha_i+\beta_j-1}{p_i\wedge q_j}, \tag{\AfrMainRes} $$
where $p_i\wedge q_j=\lim_{n\to\infty}(p_iq_j)^n$.
If $L_0$ is empty then $\Afr_0$ is a simple C$^*$--algebra.
Otherwise, for every $(i,j)\in L_0$, although
$\solim_{n\to\infty}(p_iq_j)^n=0$, there is a
unital $*$--homomorphism $\pi_{(i,j)}:\Afr_0\to\Cpx$ such that
$\pi_{(i,j)}(r_0p_i)=1=\pi_{(i,j)}(r_0q_j)$
and
$$ \Afr_{00}\eqdef\bigcap_{(i,j)\in L_0}\ker\pi_{(i,j)} $$
is simple, nonunital and with unique tracial state
$\phi(r_0)^{-1}\phi\restrict_{\Afr_{00}}$.
\endproclaim
The notation in~(\AfrMainRes) means that $\Afr=\Afr_0$ if $L_+$ is empty and
otherwise
$$ \Afr=\Afr_0\oplus\undersetbrace{|L_+|\text{ times}}\to
{\Cpx\oplus\cdots\oplus\Cpx}. $$
In addition, for each $(i,j)\in L_+$ the corresponding central summand is
$\Cpx(p_i\wedge q_j)$ and
$\tau(p_i\wedge q_j)=\alpha_i+\beta_j-1$.
The assertion that $\solim_{n\to\infty}(p_iq_j)^n=0$ when $\alpha_i+\beta_j=1$
refers to the strong--operator limit in the GNS representation of $\Afr$ with
respect to $\tau$.
Finally, $r_0=1-\sum_{(i,j)\in L_+}p_i\wedge q_j$ is the projection which is
the unit of $\Afr_0\oplus0\oplus\cdots\oplus0$.

Analogous results hold for free products of more than two finite dimensional
abelian C$^*$--algebras and for free products of direct sums of other abelian
C$^*$--algebras, (see Propositions~\FinManyFDAbel, \FinManyAzeroCdC{}
and~\InvLimIsolAtomsFP).

The following result is used in the proof of Theorem~\FPfdAbelian{} and
is also of independent interest.
\proclaim{Theorem \FPdas}
Let $A$ and $B$ be unital C$^*$--algebras with states $\phi_A$, respectively
$\phi_B$, whose GNS representations are faithful.
Let
$$ (\Afr,\phi)=(A,\phi_A)*(B,\phi_B). $$
Suppose $B\ne\Cpx$ and the centralizer of $\phi_A$ has an abelian subalgebra
$D\cong C(X)$ such that the restriction of $\phi_A$ to $D$ is given by an
atomless measure on $X$.
Then $\Afr$ is simple.
If $\phi_A$ and $\phi_B$ are traces then $\Afr$ has stable rank~1 and $\phi$ is
the unique tracial state on $\Afr$.
If one of $\phi_A$ and $\phi_B$ is not a trace then $\Afr$ has no tracial
states.
\endproclaim

\vskip3ex
\noindent{\bf\S\Preliminaries. Preliminaries.}
\vskip3ex

\proclaim{Definition \DefDiffuse}\rm
Let $A$ be a unital C$^*$--algebra and $\phi$ a state on $A$.
Let $B\cong C(X)$ be an abelian C$^*$--subalgebra of $A$ containing the unit of
$A$.
We say that $\phi$ is {\it diffuse} on $B$ if
$\phi\restrict_B$ is given by a measure
on $X$ having no atoms.
It will usually be clear from the context which state we mean, and then we will
speak simply of $B$ being a {\it diffuse abelian subalgebra} of $A$.

Given a C$^*$--algebra with state $(A,\phi)$, a {\it Haar unitary} (with
respect to $\phi$) is a unitary
element, $u\in A$, such that $\phi(u^n)=0$ for every nonzero integer $n$.
\endproclaim

\proclaim{Proposition \DASHU~\scite{\DykemaHaagerupRordam}{4.1(i)}}
Let $B$ be a unital, abelian C$^*$--algebra with state $\phi$.
Then $\phi$ is diffuse on $B$ if and only if
$B$ contains a Haar unitary (with respect to $\phi$).
\endproclaim

Recall that the {\it centralizer} of the state $\phi$ is
$\{a\in A\mid\forall x\in A\;\phi(ax)=\phi(xa)\}$.
We will often be interested in the situation when the centralizer of $\phi$
contains a Haar unitary.

An {\it ideal} of a C$^*$--algebra always means a closed, two--sided ideal.

For unital
C$^*$--algebras $A_1,A_2,\ldots,A_n$,
it is an obvious fact that
$A_1\oplus A_2\oplus\cdots\oplus A_n$
has stable rank~1 if and only if for each $j$ $A_j$ has stable rank~1.
In addition, we will make use of the following results, due to Rieffel.
\proclaim{Proposition \TensorMn~\scite{\RieffelZZsr}{3.3}}
Let $n\in\Naturals$ and let $A$ be a C$^*$--algebra.
Then $A$ has stable rank~1 if and only if $A\otimes M_n(\Cpx)$ has stable
rank~1.
\endproclaim
The following result follows from~\scite{\RieffelZZsr}{4.4 and 4.11} together
with the fact that in a finite dimensional C$^*$--algebra $B$, the left
invertible elements are invertible, hence the connected stable rank of
$B$ is one.
\proclaim{Proposition \ShortExactSeq}
Let $A$ be a C$^*$--algebra with an ideal, $J$, such that $A/J$ is finite
dimensional.
Then $A$ has stable rank~1 if and only if $J$ has stable rank~1.
\endproclaim

Recall that a hereditary C$^*$--subalgebra, $B$, of a C$^*$--algebra, $A$, is
said to be {\it full} if there is no closed, proper, two--sided ideal of $A$
containing $B$.
\proclaim{Proposition \AheredBsrone}
Let $A$ be a C$^*$--algebra with countable approximate identity.
Take $h\in A$, $h\ge0$ and let $B$ be the hereditary subalgebra,
$\overline{hAh}$, of $A$.
Suppose that $B$ is full in $A$.
Then
\roster
\item"(i)" $A$ has stable rank~1 if and only if $B$ has stable rank~1.
\item"(ii)" If $B$ has unique tracial state then $A$ has at most one tracial
state.
\endroster
\endproclaim
\demo{Proof}
For~(i), note that
$B$ has a countable approximate identity for itself because $h$ is strictly
positive in $B$, (see~\scite{\BrownLGZZStIsoHered}{p\. 327}).
Thus, by~\cite{\BrownGreenRieffel}, $A$ and $B$ are stably isomorphic, i.e\.
$A\otimes\KHil\cong B\otimes\KHil$, where $\KHil$ is the algebra of compact
operators on separable Hilbert space.
But \scite{\RieffelZZsr}{3.6} states that
$$ \sr(A)=1\qquad\Leftrightarrow\qquad\sr(A\otimes\KHil)=1 $$
and similarly for $B$, hence
$$ \sr(A)=1\qquad\Leftrightarrow\qquad\sr(B)=1. $$

To see~(ii), note that $\lspan\{xhahy\mid a,x,y\in A\}$ is dense in $A$.
If $\tau$ is a tracial state on $A$ then
$\tau(xhahy)=\tau(h^{1/2}ahyxh^{1/2})$ and $h^{1/2}ahyxh^{1/2}\in B$, so $\tau$
is determined by $\tau\restrict_B$.
\QED

We will say that a positive element, $h\in A$ is {\it full} if the hereditary
subalgebra $\overline{hAh}$ is full in $A$.
The following fact is easy to show.
\proclaim{Proposition~\FullSimple}
Let $A$ be a C$^*$--algebra and let $B$ be a full hereditary C$^*$--subalgebra
of $A$.
Then $A$ is simple if and only if $B$ is simple.
\endproclaim

In fact, (see~\cite{\RieffelZZMoritaEqOpAlg}), it is well--known that the
representation theories of a C$^*$--algebra $A$ and its hereditary
C$^*$--subalgebra $B$ are equivalent.

The reduced free product of two two--dimensional C$^*$--algebras is the most
transparent nontrivial free product one can consider.
It is understood completely and described in the proposition below.
This description is the starting point for our investigation into reduced free
products of more general finite dimensional abelian algebras.

\proclaim{Proposition \TwoProj}
Let $1>\alpha\ge\beta\ge\frac12$ and let
$$ (\Afr,\tau)=(\smdp\Cpx\alpha p\oplus\smdp\Cpx{1-\alpha}{1-p})
*(\smdp\Cpx\beta q\oplus\smdp\Cpx{1-\beta}{1-q}). $$
If $\alpha>\beta$ then
$$ \Afr=\smdp\Cpx{\alpha-\beta}{p\wedge(1-q)}\oplus
C([a,b],M_2(\Cpx))\oplus\smdp\Cpx{\alpha+\beta-1}{p\wedge q}, $$
for some $0<a<b<1$.
Furthermore, in the above picture
$$ \align
p&=1\oplus\left(\smallmatrix1&0\\0&0\endsmallmatrix\right)\oplus1 \\
q&=0\oplus\left(\smallmatrix t&\sqrt{t(1-t)}\\
\sqrt{t(1-t)}&1-t\endsmallmatrix\right)
\oplus1
\endalign $$
and the faithful trace $\tau$ is given by the indicated weights on the
projections $p\wedge(1-q)$ and $p\wedge q$, together with an atomless measure
whose support is $[a,b]$.

If $\alpha=\beta>\frac12$ then
$$ \Afr=\{f:[0,b]\to M_2(\Cpx)\mid f\text{ continuous and }f(0)
\text{ diagonal }\}\oplus\smdp\Cpx{\alpha+\beta-1}{p\wedge q}, $$
for some $0<b<1$.
Furthermore, in the above picture
$$ \align
p&=\left(\smallmatrix1&0\\0&0\endsmallmatrix\right)\oplus1 \\
q&=\left(\smallmatrix t&\sqrt{t(1-t)}\\\sqrt{t(1-t)}&1-t\endsmallmatrix\right)
\oplus1
\endalign $$
and the faithful trace $\tau$ is given by the indicated weight on the
projection $p\wedge q$, together with an atomless measure
on $[0,b]$.

If $\alpha=\beta=\frac12$ then
$$ \Afr=\{f:[0,1]\to M_2(\Cpx)\mid f\text{ continuous and }f(0)\text{ and }f(1)
\text{ diagonal }\}. $$
Furthermore, in the above picture
$$ \align
p&=\left(\smallmatrix1&0\\0&0\endsmallmatrix\right) \\
q&=\left(\smallmatrix t&\sqrt{t(1-t)}\\\sqrt{t(1-t)}&1-t\endsmallmatrix\right)
\endalign $$
and the faithful trace $\tau$ is given by an atomless measure
whose support is $[0,1]$.
\endproclaim
\demo{Proof}
Once the traces of $p$ and $q$ are known, the C$^*$--algebra $\Afr$ and the
trace $\tau$ are determined by $\tau$
composed with the functional calculus of $pqp$.
This, in turn, is computed using Voiculescu's multiplicative free
convolution~\cite{\VoiculescuZZMult}.
To see this proof in more detail, see~\cite{\DykemaZZSNU} or the
proof of the analogous result for von Neumann algebras
in~\scite{\DykemaZZFreeDim}{1.1}.
\QED

The following proposition is a variation on Theorem~1.2
of~\cite{\DykemaZZFreeDim} and is proved similarly.

\proclaim{Proposition \FreeProdDirectSum}
Let $A=A_1\oplus A_2$ be a direct sum of unital C$^*$--algebras, write
$p=1\oplus0\in A$ and let $\phi_A$ be a state on $A$, such that
$0<\alpha\eqdef\phi_A(p)<1$.
Let $B$ be a unital C$^*$--algebra with state $\phi_B$ and let
$(\Afr,\phi)=(A,\phi_A)*(B,\phi_B)$.
Let $\Afr_1$ be the C$^*$--subalgebra of $\Afr$ generated by
$(0\oplus A_2)+\Cpx p\subseteq A$ together with $B$.
We abbreviate this by writing
$$ \align
(\Afr,\phi)&=(\smdp{A_1}\alpha p\oplus\smdp{A_2}{1-\alpha}{1-p})*(B,\phi_B) \\
\cup\quad\\
(\Afr_1,\phi\restrict_{\Afr_1})&
=(\smdp\Cpx\alpha p\oplus\smdp{A_2}{1-\alpha}{1-p})*(B,\phi_B).
\endalign $$
Then $p\Afr p$ is generated by $p\Afr_1p$ and $A_1\oplus0\subseteq A$, which
are free in $(p\Afr p,\frac1\alpha\phi\restrict_{p\Afr p})$.
\endproclaim

The next elementary lemma will come in handy.
\proclaim{Lemma \notCpx}
Let $B$ be a unital C$^*$--algebra and $\phi$ a state on $B$ whose GNS
representation is faithful.
If $|\phi(u)|=1$ for every unitary, $u\in B$, then $B=\Cpx$.
\endproclaim
\demo{Proof}
Let the defining embedding $B\hookrightarrow L^2(B,\phi)$ be denoted
$b\mapsto\bh$.
Let $\Uc(B)$ denote the unitary group of $B$.
Whenever $u\in\Uc(B)$ then $\nm\uh=1$ but also $|\langle\uh,1\rangle|=1$, so
$\uh=\alpha\onehat$, some $\alpha\in\Tcirc$.
But $L^2(B,\phi)=\clspan\{\uh\mid u\in\Uc(B)\}$, so $L^2(B,\phi)$ is
one--dimensional.
This implies $B=\Cpx$.
\QED

\vskip3ex
\noindent{\bf\S\InPresOfDiff. When in the presence of one spread thin.}
\nopagebreak
\vskip3ex
\nopagebreak
Let $A$ and $B$ be unital C$^*$--algebras with states, $\phi_A$ and $\phi_B$,
respectively, whose GNS representations are faithful and let
$$ (\Afr,\phi)=(A,\phi_A)*(B,\phi_B). $$
In this section we will prove that if
the centralizer of $\phi_A$ contains a diffuse abelian subalgebra
then $\Afr$ is simple and if, furthermore, $\phi_A$ and $\phi_B$ are traces
then $\Afr$ has stable rank~1.
The diffuse abelian subalgebra is, if you like, ``one spread thin.''
By~\scite{\DykemaHaagerupRordam}{4.1(i)}, (see Proposition~\DASHU{} above) this
condition is equivalent to the
centralizer of $\phi_A$ containing a Haar unitary, $u$.

Denote by $\Afr_0$ the norm dense $*$--subalgebra of $\Afr$ that is generated
by $A\cup B$.
Then using the standard notation $\Ao=\ker\phi_A$ and $\Bo=\ker\phi_B$,
every element, $x$, of $\Afr_0$ can be written
$x=x_0+x_1$ where $x_0\in A$ and
$$ x_1=\sum_{j=1}^Na_o^{(j)}b_1^{(j)}a_1^{(j)}b_2^{(j)}a_2^{(j)}\cdots
b_{n(j)}^{(j)}a_{n(j)}^{(j)} \tag{\xoneab} $$
with $N\in\Naturals$, $n(j)\in\Naturals$, $a_0^{(j)},a_{n(j)}^{(j)}\in A$,
$a_1^{(j)},\ldots,a_{n(j)-1}^{(j)}\in\Ao$,
$b_1^{(j)},\ldots,b_{n(j)}^{(j)}\in\Bo$.
Expressed another way,
$$ x_1\in\lspan\left(\bigcup_{n=1}^\infty A\undersetbrace{n\text{ times }\Bo}
\to{\Bo\Ao\cdots\Bo\Ao\Bo}A\right). $$

We begin with a technical lemma.
Let $u$ be a Haar unitary in the centralizer of $\phi\restrict_A$ and write
$$ \align
u^{-\Naturals}&\eqdef\{u^{-k}\mid k\in\Naturals\} \\
u^\Naturals&\eqdef\{u^k\mid k\in\Naturals\}.
\endalign $$

\proclaim{Lemma \TechnicalLemma}
With notation as above, suppose that $B\ne\Cpx$ and the centralizer of $\phi_A$
contains a Haar unitary, $u$.
Given $\epsilon>0$ and $x\in\Afr$ such that $\phi(x)=0$, there is a
unitary, $z\in\Afr$ such that $z^*xz$ differs in norm by no more than
$\epsilon$ from a finite linear combination of elements of
$$ \Theta\eqdef\bigcup_{n=1}^\infty
u^{-\Naturals}\undersetbrace{n\text{ times }\Bo}
\to{\Bo\Ao\cdots\Bo\Ao\Bo}u^\Naturals. \tag{\uNBAB} $$
\endproclaim
\demo{Proof}
Since $\Afr_0$ is dense in $\Afr$, we may assume without loss of generality
that $x\in\Afr_0$.
By Lemma~\notCpx{}
there is a unitary element,
$v\in B$, such that $0\le\phi_B(v)<1$.
Let $c_0=\phi_B(v)$, $c_1=\sqrt{1-c_0^2}$ and $y=(v-c_01)/c_1$, so that $1$ and
$y$ are orthonormal in $L^2(B,\phi_B)$.
Let $n,k\in\Naturals$ and let $z=(u^kv)^nu^k$.
Write $x=x_0+x_1$ with $x_0\in A$ and $x_1$ as in~(\xoneab).

We first concern ourselves with $z^*x_1z$.
Writing $x_1$ as in~(\xoneab), let $\eta>0$.
Since $(u^p)_{p\in\Integers}$ is an orthonormal family in
$L^2(A,\phi_A)$, we have
$$ \forall a\in A\qquad\lim_{p\to\infty}\phi_A(au^p)=0
=\lim_{p\to\infty}\phi_A(au^{-p}). \tag{\limLtwo} $$
Using~(\limLtwo) we see that if $k$ is large enough then for every positive
integer $p$ and every $j$,
$|\phi_A(u^{-pk}a_0^{(j)})|<\eta$ and
$|\phi_A(a_{n(j)}^{(j)}u^{pk})|<\eta$.
Since $v=c_01+c_1y$, we have
$$ a_{n(j)}^{(j)}z=\sum_{\delta_1,\ldots,\delta_n\in\{0,1\}}
c_{\delta_1}\cdots c_{\delta_n}a_{n(j)}^{(j)}
u^ky^{\delta_1}\cdots u^ky^{\delta_n}u^k, $$
where
$$ \align
y^{-\delta}&=\cases y^*&\text{if }\delta=1\\1&\text{if }\delta=0\endcases \\
y^\delta&=\cases y&\text{if }\delta=1\\1&\text{if }\delta=0.\endcases
\endalign $$
If not all the $\delta_j$ are zero than
$a_{n(j)}^{(j)}u^ky^{\delta_1}\cdots u^ky^{\delta_n}u^k$ differs in norm by at
most $\eta\nm{y}^m$ (where $m$ is the number of $1\le j\le n$ for which
$\delta_j=1$) from an element of
$$ \undersetbrace{m\text{ times }\Ao y}\to
{\Ao y\Ao y\cdots\Ao y}u^\Naturals. $$
Since $\nm y<(1+c_0)/c_1$, we obtain that $a_{n(j)}^{(j)}z$ differs in norm by
at most $c_0^n\nm{a_{n(j)}^{(j)}}+(1+2c_0)^n\eta$ from an element of
$$ \lspan\left(\bigcup_{m=1}^n\undersetbrace{m\text{ times }\Ao y}\to
{\Ao y\Ao y\cdots\Ao y}u^\Naturals\right). $$
Similarly, $z^*a_0^{(j)}$ differs in norm by at most
$c_0^n\nm{a_0^{(j)}}+(1+2c_0)^n\eta$ from an element of
$$ \lspan\left(\bigcup_{m=1}^nu^{-\Naturals}
\undersetbrace{m\text{ times }y^*\Ao}\to{y^*\Ao y^*\Ao\cdots y^*\Ao}\right). $$
Therefore $z^*x_1z$ dffers in norm by no more than
$$ \sum_{j=1}^N\bigl(c_0^n\nm{a_0^{(j)}}+(1+2c_0)^n\eta\bigr)\;
\nm{b_1^{(j)}a_1^{(j)}b_2^{(j)}a_2^{(j)}\cdots b_{n(j)}^{(j)}}\;
\bigl(c_0^n\nm{a_{n(j)}^{(j)}}+(1+2c_0)^n\eta\bigr) $$
from a finite
linear combination of elements from $\Theta$.
Thus, if $n$ is chosen large enough and then $k$ is chosen large enough, then
$z^*x_1z$ can be made arbitrarily close to a finite
linear combination of elements from $\Theta$.

We now examine $z^*x_0z$.
Using again $v=c_01+c_1y$ we have
$$ z^*x_0z=\sum_{\delta_1,\ldots,\delta_{2n}\in\{0,1\}}
c_{\delta_1}\cdots c_{\delta_{2n}}
u^{-k}y^{-\delta_1}u^{-k}\cdots y^{-\delta_n}u^{-k}x_0
u^ky^{\delta_{n+1}}\cdots u^ky^{\delta_{2n}}u^k. \tag{\zcongxzero} $$
We first concentrate on the $2^{n+1}-1$ terms
when either $\delta_1=\delta_2=\cdots=\delta_n=0$ or
$\delta_{n+1}=\delta_{n+2}=\cdots=\delta_{2n}=0$.
The sum over these terms is equal to
$$ c_0^n(u^{-(n+1)k}x_0z+z^*x_0u^{(n+1)k}-c_0^nu^{-(n+1)k}x_0u^{(n+1)k}), $$
which has norm no greater than $c_0^n\nm{x_0}(2+c_0^n)$.
This can be made
arbitrarily small by choosing $n$ large enough (independently of $k$).
Each of the remaining $2^{2n}-2^{n+1}+1$ terms of~(\zcongxzero) is of the form
$$ c_0^lc_1^{l'}
u^{-r_pk}y^*\cdots u^{-r_2k}y^*u^{-r_1k}y^*u^{-r_0k}x_0
u^{s_0k}yu^{s_1k}yu^{s_2k}\cdots yu^{s_qk}, \tag{\usualterms} $$
where $l',p,q,r_j,s_j$ are positive integers and $l\ge0$.
Clearly
$$ \phi(u^{-r_0k}x_0u^{s_0k})=\phi(x_0u^{(s_0-r_0)k}). \tag{\urxus} $$
If $r_0=s_0$ then $\phi(u^{-r_0k}x_0u^{s_0k})=\phi(x_0)=0$ and hence the
term~(\usualterms) is an element of $\Theta$.
Using~(\limLtwo) we see that by choosing $k$ large enough, each
quantity~(\urxus) can be made
arbitrarily small and hence each of the terms~(\usualterms) can be made
arbitrarily close to an element of $\Theta$.
Thus if $n$ is chosen large enough and then $k$ is chosen large enough, then
$z^*x_0z$ can be made arbitrarily close to a finite linear combination of
elements from $\Theta$.

Considering the above analyses for $z^*x_1z$ and $z^*x_0z$ at the same time, we
can choose $n$ large enough and then $k$ large enough so that $z^*xz$ is
arbitrarily close to a finite linear combination of elements from $\Theta$.
\QED

\proclaim{Proposition \FPHaarUnitaryDixmier}
Let $A$ and $B$ be unital C$^*$--algebras with states, $\phi_A$ and $\phi_B$,
respectively, whose GNS representations are faithful.
Let
$$ (\Afr,\phi)=(A,\phi_A)*(B,\phi_B). $$
Suppose the centralizer of $\phi_A$ has diffuse abelian subalgebra and
$B\ne\Cpx$.
Then for every
$x\in\Afr$ and
$\epsilon>0$ there are $n\in\Naturals$ and unitaries
$z_1,\ldots,z_n\in\Afr$ such that
$$ \nm{\phi(x)1-\frac1n\sum_{r=1}^nz_r^*xz_r}<\epsilon. \tag{\dixmierprop} $$
Consequently, $\Afr$ is simple.
Moreover, if both $\phi_A$ and $\phi_B$ are traces then $\phi$ is the unique
tracial state on $\Afr$.
If one or both of $\phi_A$ and $\phi_B$ is not a trace then $\Afr$ has no
tracial states.
\endproclaim
\demo{Proof}
To prove the existence of $z_r$ such that~(\dixmierprop) holds,
we may without loss of generality assume that $x=x^*$ and $\phi(x)=0$,
and we may replace $x$ by a unitary conjugate of itself.
Let $u$ be a Haar unitary in the centralizer of $\phi_A$
Employing Lemma~\TechnicalLemma, we may assume that
$x\in\lspan\Theta$, (see~(\uNBAB)).
Now we will find
$z_1,\ldots,z_5\in u^\Naturals$ such that
$$ \nm{\frac15\sum_{r=1}^5z_r^*xz_r}\le\frac{49}{50}\nm x.
\tag{\DixmierFive} $$
These will be found using the technique of~\cite{\PowersZZCsFtwo}.
With notation similar to~(\xoneab), we have
$$ x=\sum_{j=1}^Nu^{-l_j}
b_1^{(j)}a_1^{(j)}\cdots b_{n(j)-1}^{(j)}a_{n(j)-1}^{(j)}b_{n(j)}^{(j)}u^{m_j},
$$
for $l_j,m_j\in\Naturals$.
Let $K=\max(\bigcup_{j=1}^N\{l_j,m_j\})+1$.
For $1\le r\le5$, let $z_r=u^{rK}$.
Let $\Mfr_r$ be the closed subspace of $L^2(\Afr
,\phi)$ spanned by all words of
the form $u^{-k}b_1a_1\cdots b_na_n$ with $k\in\Naturals$, $(r-1)K<k\le rK$,
$n\in\Naturals\cup\{0\}$, $b_1,\ldots,b_n\in\Bo$, $a_1,\ldots,a_{n-1}\in\Ao$
and $a_n\in A$.
Clearly $p\ne q$ implies $\Mfr_p\perp\Mfr_q$.
Since $z_r^*xz_r$ is a finite sum of words whose left--most letter lies in
$\{u^{-k}\mid(r-1)K<k\le rK\}$ and whose right--most letter lies in
$\{u^k\mid(r-1)K<k\le rK\}$, we have
$$ z^*_rxz_r(\Mfr_r^\perp)\subseteq\Mfr_r. $$
Denote by $E_r$ the projection from $L^2(\Afr,\phi)$ onto
$\Mfr_r$.
Given a unit vector, $\xi\in L^2(\Afr,\phi)$, there is some $1\le p\le5$ for
which $\nm{E_p\xi}^2\le\frac15$.
Thus
$$ |\langle\frac15\sum_{r=1}^5z^*_rxz_r\xi,\xi\rangle|\le\frac45\nm x
+\frac15|\langle z_p^*xz_p\xi,\xi\rangle| $$
and since $(1-E_p)z_p^*xz_p(1-E_p)=0$ we have
$$ |\langle z_p^*xz_p\xi,\xi\rangle|=|\langle z_p^*xz_pE_p\xi,\xi\rangle|+
|\langle z_p^*xz_p(1-E_p)\xi,E_p\xi\rangle|
\le2\nm x\,\nm{E_p\xi}\le\frac2{\sqrt 5}\nm x. $$
Hence
$$
|\langle\frac15\sum_{r=1}^5z^*_rxz_r\xi,\xi\rangle|\le
\left(\frac45+\frac2{\sqrt 5}\right)\nm x<\frac{49}{50}\nm x. $$
This implies~(\DixmierFive).

To finish the proof of~(\dixmierprop), note that the
element
$\frac15\sum_{r=1}^5z^*_rxz_r$ obtained above is again in $\lspan\Theta$.
Hence, by repeating this process as many times as necessary, for any
$\epsilon>0$ there are $n\in\Naturals$ and
$z_1,\ldots,z_n\in u^\Naturals$ such that
$\nm{\frac1n\sum_{r=1}^nz_r^*xz_r}<\epsilon$.

Now the remaining facts follow by standard arguments
of \cite{\PowersZZCsFtwo} and~\cite{\Avitzour}.
Indeed, suppose $\Ic$ is a nonzero, two--sided, closed ideal of $\Afr$.
Let $a\in\Ic\backslash\{0\}$.
Since the GNS representation of $\Afr$ associated to $\phi$ is faithful, there
must be $b\in\Afr$ such that $\phi(b^*a^*ab)\ne0$.
Then from~(\dixmierprop), it follows that $\phi(b^*a^*ab)1\in\Ic$, hence
$\Ic=\Afr$ and consequently $\Afr$ is simple.

The property described at~(\dixmierprop) implies that any tracial state on
$\Afr$ must be equal to $\phi$.
If both $\phi_A$ and $\phi_B$ are traces, then the free product state, $\phi$
is also a trace, and is thus the unique tracial state.
If one of $\phi_A$ and $\phi_B$ is not a trace, then neither is $\phi$ a trace,
hence $\Afr$ has no tracial states.
\QED

\proclaim{Lemma \TechnicalLemmatwo}
Let $x\in\Afr$ and let $\epsilon>0$.
Then there are unitaries, $z_1,z_2\in\Afr$ such that $\nm{z_1xz_2-x'}<\epsilon$
for some $x'\in\lspan\Theta$, with $\Theta$ as in~(\uNBAB).
\endproclaim
\demo{Proof}
Let $u\in A$ be a Haar unitary and $v\in B$ a unitary such that
$0\le\phi_B(v)<1$.
By Lemma~\TechnicalLemma{} there is a unitary, $z\in\Afr$ such that
$$ \nm{z^*xz-(\phi(x)1+x'')}<\epsilon/2, \tag{\zxz} $$
where $x''\in\lspan\Theta$.
Writing $v=c_01+c_1y$ like in the proof of Lemma~\TechnicalLemma, we see that
$(u^*v)^pu^p$ differs in norm by no more that $c_0^p$ from an element of
$\lspan\Theta$.
We similarly see that $(u^*v)^px''\in\lspan\Theta$.
Let $p$ be so large that $c_0^p<\epsilon/(2+2|\phi(x)|)$.
Let $z_1=(u^*v)^pz^*$ and $z_2=zu^p$.
Then from~(\zxz) we have
$$ \nm{z_1xz_2-\bigl(\phi(x)(u^*v)^pu^p+(u^*v)^px''u^p\bigr)}
<\epsilon/2, $$
and from the above discussion there is $x'\in\lspan\Theta$ such that
$$ \nm{\phi(x)(u^*v)^pu^p+(u^*v)^px''u^p-x'}<\epsilon/2. $$
Hence $\nm{z_1xz_2-x'}<\epsilon$.
\QED

The proof of the following proposition uses ideas
from~\cite{\DykemaHaagerupRordam}.
\proclaim{Proposition \FPHaarUnitarySRone}
Let $A$ and $B$ be unital C$^*$--algebras with faithful, tracial states,
$\tau_A$ and $\tau_B$, respectively.
Let
$$ (\Afr,\tau)=(A,\tau_A)*(B,\tau_B). $$
Suppose $A$ has a diffuse abelian subalgebra and $B\ne\Cpx$.
Then
$\Afr$ has stable rank $1$, i.e\. the set of invertible elements of
$\Afr$ is dense in $\Afr$.
\endproclaim
\demo{Proof}
Suppose for contradiction that the set of invertibles in $\Afr$, denoted
$\GL(\Afr)$, is not dense in $\GL(\Afr)$.
Then, by~\scite{\RordamZZUnitaryRank}{2.6}, there is $x\in\Afr$ such that
$\nm x=1$ and $\dist(x,\GL(\Afr))=1$.
We must have $\nm x_2<1$ since $\nm x_2=1$ would imply that $x$ is unitary.
Let $\epsilon=1-\nm x_2$.
Let $u$ be a Haar unitary in $A$ and let $v\in B$ be a unitary such that
$0\le\tau_B(v)<1$.
By Lemma~\TechnicalLemmatwo{} there are $n\in\Naturals$ and unitaries,
$z_1,z_2\in\Afr$, such that
$\nm{z_1xz_2-x'}<\epsilon/8$ for some $x'\in\lspan\Psi_n$,
where 
$$ \Psi_n=\{u^lb_1a_1\cdots b_{k-1}a_{k-1}b_ku^m
\mid k,l,m\in\Naturals,\,l,m<n,\,b_j\in\Bo,a_j\in\Ao\}. $$
Let $p$ be so large that $c_0^p<\epsilon/(8(\nm x+\epsilon))$.
By writing $v=c_01+c_1y$ as in the proof of Lemma~\TechnicalLemma, we see that
$$ \nm{(u^nv)^pu^nx'-x''}<\epsilon/8 $$
for some $x''\in\lspan\Psi_{n,p}$ where
$$ \Psi_{n,p}=\{u^{nl}b_1a_1\cdots b_{k-1}a_{k-1}b_ku^m
\mid k,l,m\in\Naturals,\,l<p,\,m<n,\,b_j\in\Bo,\,a_j\in\Ao\}. $$
For $q\in\Naturals$ let $E_q:A\to A$ be
$$ E_q(a)=\sum_{j=q+1}^{q+np}u^j\langle\ah,(u^j)\hat{\;}\rangle, $$
where we denote the defining embedding
$A\hookrightarrow L^2(A,\tau_A)$ by $a\mapsto\ah$.
Since
$\lim_{j\to\infty}\langle\ah,(u^j)\hat{\;}\rangle=0$
we have $\lim_{q\to\infty}\nm{E_q(a)}=0$ for every $a\in A$.
Therefore, there is $q$, a positive multiple of $n$, such that
$\nm{x^{(3)}-x''}<\epsilon/8$ for some $x^{(3)}$ in $\lspan\Psi_{n,p,q}'$,
where
$$ \Psi_{n,p,q}'=\{u^{nl}b_1a_1\cdots b_{k-1}a_{k-1}b_ku^m
\mid k,l,m\in\Naturals,\,l<p,\,m<n,\,b_j\in\Bo,\,a_j\in\Ao,E_q(a_j)=0\}. $$
Let $X_A$ be a standard orthonormal basis
(see~\scite{\DykemaHaagerupRordam}{\S2}) for $(A,\tau_A)$ containing
$\{u,u^2,u^3,\ldots,u^{q+np}\}$ and let
$X_B$ be a standard orthonormal basis for $(B,\tau_B)$.
Let $Y=X_A*X_B$ be the resulting
free product standard orthonormal basis for $(\Afr,\tau)$.
(Note that, by definition, $Y\backslash\{1\}$ is the set of all reduced words
in $X_A\backslash\{1\}$ and $X_B\backslash\{1\}$.)
Then there is $x^{(4)}\in\lspan Y_{n,p,q}'$ such that
$\nm{u^qx^{(3)}-x^{(4)}}<\epsilon/8$,
where $Y_{n,p,q}'$ is the subset of $Y$ defined by
$$ \align
Y_{n,p,q}'=\bigg\{u^{q+nl}b_1a_1\cdots b_{k-1}a_{k-1}b_ku^m
\bigg|\,&l,m,k\in\Naturals,\,l<p,\,m<n,\, \\
&b_j\in X_B\backslash\{1\}, \\
&a_j\in X_A\backslash\{1,u^{q+1},u^{q+2},\ldots,u^{q+np}\}\bigg\}.
\endalign $$
Now we see that, since
no cancellation occurs when we multiply elements of $Y'_{n,p,q}$, (but only
``$u$ on $u$ contact''), whenever $w_1,\ldots,w_m\in Y'_{n,p}$ we have
$$ \nm{w_1w_2\cdots w_m}_2
=\nm{w_1}_2\nm{w_2}_2\cdots\nm{w_m}_2. $$
Moreover, if $w_1w_2\cdots w_m=w_1'w_2'\cdots w_m'$ for
any $m\in\Naturals$ and
$w_j,w'_j\in Y_{n,p,q}'$, then
$w_1'=w_1$, $w_2'=w_2$, $\ldots$, $w_m'=w_m$.
The reason for this is that when we take the reduced word of
$w_1w_2\ldots w_m$, a letter $u^k$ for $q\le k\le q+np$ appears at
every boundary where $w_j$ touches $w_{j+1}$, $(1\le j\le m-1)$, and
nowhere else,
and writing $u^k=u^ru^{q+ln}$ for $l,r\in\Naturals$ and $r\le n$, we see that
$u^r$ must have been the last letter in $w_j$ and $u^{q+ln}$ must have been the
first letter in $w_{j+1}$, so we can recover the list of letters,
$w_1,w_2,\ldots,w_m$ from their product.
Thus we see that
$\nm{(x^{(4)})^m}_2=\nm{x^{(4)}}_2^m$ for every $m\in\Naturals$ and
(see~\scite{\DykemaHaagerupRordam}{3.2}) $K((x^{(4)})^m)=K(x^{(4)})$.
Now we argue as in the proof of~\scite{\DykemaHaagerupRordam}{3.8} to show that
the spectral radius of $x^{(4)}$, denoted $r(x^{(4)})$, is no greater than
$\nm{x^{(4)}}_2$.
Indeed, let $q$ be the largest block length of the words in the support of
$x^{(4)}$, so that, in the notation of~\scite{\DykemaHaagerupRordam}{2.2},
$$ x^{(4)}\in\lspan\bigcup_{j=1}^q Y_j. $$
Then, by~\scite{\DykemaHaagerupRordam}{3.5},
$$ \forall m\in\Naturals\qquad\nm{(x^{(4)})^m}
\le(2mk+1)^{3/2}K(x^{(4)})\nm{x^{(4)}}_2^m, $$
where $K(x^{(4)})$ is a constant.
Hence
$$ r(x^{(4)})=\liminf_{m\to\infty}\nm{(x^{(4)})^m}^{1/m}\le\nm{x^{(4)}}_2. $$
Therefore $\dist(x^{(4)},\GL(\Afr))\le\nm{x^{(4)}}_2$.
But $\nm{x^{(4)}-u^k(u^nv)^pu^nz_1xz_2}<\epsilon/2$, so
$$ \align
\dist(x,\GL(\Afr))&=\dist(u^k(u^nv)^pu^nz_1xz_2,\GL(\Afr))
\le\nm{x^{(4)}-u^k(u^nv)^pu^nz_1xz_2}+\nm{x^{(4)}}_2 \\
&<\epsilon/2+\nm{x^{(4)}-u^k(u^nv)^pu^nz_1xz_2}_2+\nm{u^k(u^nv)^pu^nz_1xz_2}_2
<\epsilon+\nm x_2=1,
\endalign $$
contradicting the choice of $x$.
\QED

Note that Propositions~\FPHaarUnitaryDixmier{} and~\FPHaarUnitarySRone{}
combine to prove Theorem~\FPdas.

\proclaim{Proposition \ExistsDAS}
Let $0<\alpha<1$, let $A$ be a unital C$^*$--algebra with state $\phi_A$ whose
GNS representation is faithful and let
$$ (\Afr,\phi)=(\smdp\Cpx\alpha p\oplus\smdp\Cpx{1-\alpha}{1-p})*(A,\phi_A). $$
Suppose the centralizer of $\phi_A$ has a diffuse abelian subalgebra.
Then the centralizer of $\phi\restrict_{p\Afr p}$ has a diffuse abelian
subalgebra.
\endproclaim
\demo{Proof}
Let $u$ be a Haar unitary in the centralizer of $\phi_A$, let $q=u^*pu$ and let
$B$ be the C$^*$--algebra generated by $\{1,p,q\}$.
Then $p$ and $q$ are free, so
$$ (B,\phi\restrict_B)=(\smdp\Cpx\alpha p\oplus\smdp\Cpx{1-\alpha}{1-p})
*(\smdp\Cpx\alpha q\oplus\smdp\Cpx{1-\alpha}{1-q}). $$

\proclaim{Case I(\ExistsDAS)} $\alpha<1/2$. \endproclaim
Then by Proposition~\TwoProj{} we have for some $0<b<1$ that
$$ B\cong\{f:[0,b]\to M_2(\Cpx)\mid f\text{ continuous and }
f(0)\text{ diagonal }\}\oplus\smdp\Cpx{1-2\alpha}{(1-p)\wedge(1-q)}, $$
$pBp\cong C([0,b])$ and $\phi\restrict_{pBp}$ is given by an atomless measure
on $[0,b]$.
Thus $pBp$ is a diffuse abelian subalgebra of the centralizer of
$\phi\restrict{pAp}$.

\proclaim{Case II(\ExistsDAS)} $\alpha=1/2$. \endproclaim
This is just as in Case I, except now
$$ B\cong\{f:[0,1]\to M_2(\Cpx)\mid f\text{ continuous and }
f(0)\text{ and }f(1)\text{ diagonal }\}. $$

\proclaim{Case III$_{\bold n}$(\ExistsDAS)} $(n\in\Naturals)$.
$1-2^{-(n-1)}<\alpha\le1-2^{-n}$.
\endproclaim
By induction on $n$.
The case III$_1$ reduces to Cases~I and~II.
Let $n>1$.
Then $\alpha>1/2$ and by Proposition~\TwoProj{} there is some $0<b<1$ such that
$$ B\cong\{f:[0,b]\to M_2(\Cpx)\mid f\text{ continuous and }
f(0)\text{ diagonal }\}\oplus\smdp\Cpx{2\alpha-1}{p\wedge q}, $$
$(p-p\wedge q)B(p-p\wedge q)\cong C([0,b])$ and the restriction of $\phi$ to
$(p-p\wedge q)B(p-p\wedge q)$ is given by an atomless measure on $[0,b]$.
Hence it will suffice to find a diffuse abelian subalgebra of the centralizer
of $\phi\restrict_{(p\wedge q)\Afr(p\wedge q)}$, because adding it to
$(p-p\wedge q)B(p-p\wedge q)$ will give a diffuse abelian subalgebra of the
centralizer of $\phi\restrict_{p\Afr p}$.

We claim that the family $\{p,q,u^2\}$ is free in $(\Afr,\phi)$.
Indeed, it suffices to show that every reduced word in $p-\phi(p)1$,
$q-\phi(q)1$ and nonzero powers of $u^2$ evaluates to zero under $\phi$.
However, rewriting each $q-\phi(q)1$ as $u^*(p-\phi(p)1)u$, we see that each
such word is equal to a word in $p-\phi(p)1$ and nonzero powers of $u$.
From the freeness of $p$ and $u$, it follows that this word evaluates to zero
under $\phi$.

Hence $p\wedge q$ and $u^2$ are free.
Letting $D$ be the C$^*$--algebra generated by $\{p\wedge q,u^2\}$, we have
$$ (D,\phi\restrict_D)\cong(\smdp\Cpx{2\alpha-1}{p\wedge q}\oplus\Cpx)
*(C^*(\Integers),\tau\restrict_\Integers). $$
Since $2\alpha-1\le1-2^{-(n-1)}$, by inductive hypothesis there is a diffuse
abelian subalgebra of $(p\wedge q)D(p\wedge q)$.
As remarked above, this finishes the proof.
\QED

\proclaim{Corollary \ExistsDAScorr}
Let $A$ be a C$^*$--algebra with state $\phi_A$ whose GNS
representation is faithful and let $n\in\Naturals$, $n\ge2$ and
$$
(\Afr,\phi)=(\smdp\Cpx{\alpha_1}{p_1}\oplus\smdp\Cpx{\alpha_2}{p_2}\oplus
\cdots\oplus\smdp\Cpx{\alpha_n}{p_n})*(A,\phi_A). $$
Suppose the centralizer of $\phi_A$ has a diffuse abelian subalgebra.
Then the centralizer of $\phi$ has a diffuse abelian
subalgebra containing $\{p_1,p_2,\ldots,p_n\}$.
\endproclaim
\demo{Proof}
For each $j$, using Proposition~\ExistsDAS{} and considering the subalgebra of
$\Afr$ generated by $A\cup\{p_j\}$,
we see that the centralizer of $\phi\restrict_{p_j\Afr p_j}$ has a diffuse
abelian subalgebra, $D_j$.
Then $D_1+D_2+\cdots+D_n$ is the required diffuse abelian subalgebra of the
centralizer of $\phi$.
\QED

\vskip3ex
\noindent{\bf\S\FPFinManyS.  Finite dimensional abelian algebras.}
\nopagebreak
\vskip3ex
\nopagebreak
In this section, we examine the reduced free product of (finitely many) finite
dimensional abelian C$^*$--algebras.
The methods used are reminiscent of~\cite{\DykemaZZFreeDim}.

Some words about notation are in order.
The natural notation
$$ (A,\tau_A)=\smdp\Cpx{\alpha_1}{p_1}\oplus\smdp\Cpx{\alpha_2}{p_2}
\oplus\cdots\oplus\smdp\Cpx{\alpha_n}{p_n} $$
for a finite dimensional abelian C$^*$--algebra and a faithful state was
explained just before Theorem~\FPfdAbelian.
Similarly, the notation
$$ \Afr=\smp{\Afr_0}{r_0}\bigoplus_{k}\smdp\Cpx{\gamma_k}{r_k} $$
was explained after that theorem.
Analogously, we will often write expressions like
$$ (\Afr,\phi)=(A_0\oplus\smdp\Cpx{\alpha_1}{p_1}\oplus\cdots\oplus
\smdp\Cpx{\alpha_n}{p_n})*(B,\phi_B). \tag{\NotationAzeroCdC} $$
This will mean $(\Afr,\phi)=(A,\phi_A)*(B,\phi_B)$ where
$$ A=A_0\oplus\undersetbrace{n\text{ times}}\to{\Cpx\oplus\cdots\oplus\Cpx}, $$
where $A_0$ is some C$^*$--algebra, where
$$ p_k=\undersetbrace{k\text{ times}}\to{0\oplus\cdots\oplus0}
\oplus1\oplus\undersetbrace{n-k\text{ times}}\to{0\oplus\cdots\oplus0} $$
and where the state, $\phi_A$, satisfies $\phi_A(p_k)=\alpha_k$.
In the case of~(\NotationAzeroCdC)
we will always assume that every $\alpha_k>0$, that $\sum_1^n\alpha_k<1$ and
that the GNS representation of the restriction of $\phi_A$ to
$A_0\oplus0\oplus\cdots\oplus0$ is faithful.
Usually, we will also desire that the centralizer of the restriction of
$\phi_A$ to $A_0\oplus0\oplus\cdots\oplus0$ have an abelian subalgebra on which
it is diffuse, (see Definition~\DefDiffuse).
This is conveniently expressed by writing ``the centralizer of
$\phi\restrict_{A_0}$ has a diffuse abelian subalgebra.''

\proclaim{Lemma \AzeroCCC}
Let
$$ (\Afr,\phi)=(A_0\oplus\smdp\Cpx\alpha p)*
(\smdp\Cpx{\beta_1}{q_1}\oplus\smdp\Cpx{\beta_2}{q_2}), $$
where the centralizer of $\phi\restrict_{A_0}$ has a diffuse abelian
subalgebra.
Take $\beta_1\ge\beta_2$.
Then
$$ \Afr=\cases\smp{\Afr_0}{r_0}&\text{if }\alpha+\beta_1\le1\\
\smp{\Afr_0}{r_0}\oplus\smdp\Cpx{\alpha+\beta_1-1}{p\wedge q_1}
&\text{if }\alpha+\beta_1>1,\,\alpha+\beta_2\le1\\
\smp{\Afr_0}{r_0}\oplus\smdp\Cpx{\alpha+\beta_1-1}{p\wedge q_1}
\oplus\smdp\Cpx{\alpha+\beta_2-1}{p\wedge q_2}
&\text{if }\alpha+\beta_2>1,
\endcases \tag{\AzeroCCCcases} $$
where the centralizer of $\phi\restrict_{\Afr_0}$ has a diffuse abelian
subalgebra which contains $r_0p$ and a diffuse abelian subalgebra which
contains $r_0q_1$, and where
$r_0p$ is full in $\Afr_0$.

If $\phi\restrict_{A_0}$ is a trace then the stable rank of $\Afr$ is~1.

If $\alpha+\beta_1\ne1$ and $\alpha+\beta_2\ne1$ then $\Afr_0$ is simple.
If, in addition, $\phi\restrict_{A_0}$ is a trace then
$\phi(r_0)^{-1}\phi\restrict_{\Afr_0}$
is the unique tracial state
on $\Afr_0$ and if $\phi\restrict_{A_0}$ is not a trace then $\Afr_0$ has no
tracial states.

Whenever $\alpha+\beta_i=1$ for $i\in\{1,2\}$, there is a $*$--homomorphism,
$\pi_i:\Afr_0\to\Cpx$, such that $\pi_i(r_0p)=1=\pi_i(q_i)$.

If $\alpha+\beta_1=1$ and $\alpha+\beta_2<1$ then $q_1$ is full in
$\Afr_0$ and $\ker\pi_1$ is simple.
If, in addition, $\phi\restrict_{A_0}$ is a trace then
$\phi(r_0)^{-1}\phi\restrict_{\ker\pi_1}$
is the unique tracial state
on $\ker\pi_1$ and if $\phi\restrict_{A_0}$ is not a trace then $\ker\pi_1$ has
no tracial states.

If $\alpha+\beta_1>1$ and $\alpha+\beta_2=1$ them $q_2$ is full in
$\Afr_0$ and $\ker\pi_2$ is simple.
If, in addition, $\phi\restrict_{A_0}$ is a trace then
$\phi(r_0)^{-1}\phi\restrict_{\ker\pi_2}$
is the unique tracial state
on $\ker\pi_2$ and if $\phi\restrict_{A_0}$ is not a trace then $\ker\pi_2$ has
no tracial states.

If $\alpha+\beta_1=1$ and $\alpha+\beta_2=1$ (which implies
$\alpha=\beta_1=\frac12$) then $q_1$ is full in $\ker\pi_2$ and $q_2$ is full
in $\ker\pi_1$ and $(\ker\pi_1)\cap(\ker\pi_2)$ is simple.
If, in addition, $\phi\restrict_{A_0}$ is a trace then
$\phi(r_0)^{-1}\phi\restrict_{\ker\pi_1\cap\ker\pi_2}$
is the unique tracial state
on $\ker\pi_1\cap\ker\pi_2$ and if $\phi\restrict_{A_0}$ is not a trace then
$\ker\pi_1\cap\ker\pi_2$ has no tracial states.
\endproclaim
\demo{Proof}
Let $\Afr_1$ be the C$^*$--subalgebra of $\Afr$ generated by $\{1,p,q\}$, so
$$ (\Afr_1,\phi\restrict_{\Afr_1})=
(\smdp\Cpx{1-\alpha}{1-p}\oplus\smdp\Cpx\alpha p)
*(\smdp\Cpx{\beta_1}{q_1}\oplus\smdp\Cpx{\beta_2}{q_2}). $$
By Proposition~\FreeProdDirectSum, $(1-p)\Afr(1-p)$ is isomorphic to the free
product of $(1-p)\Afr_1(1-p)$ and $A_0$.
We use Proposition~\TwoProj{} to find $\Afr_1$.
We will also use that $\Afr$ is generated by
$$ (1-p)\Afr(1-p)\cup(1-p)\Afr_1p\cup p\Afr_1p. $$

\proclaim{Case I(\AzeroCCC)} $\alpha>\beta_1$. \endproclaim
Then
$$ \Afr_1=\smdp\Cpx{\alpha+\beta_2-1}{p\wedge q_2}\oplus
\Bigl(C([a,b])\otimes M_2(\Cpx)\Bigr)
\oplus\smdp\Cpx{\alpha+\beta_1-1}{p\wedge q_1} $$
so $(1-p)\Afr(1-p)\cong C([a,b])*A_0$ is simple by
Proposition~\FPHaarUnitaryDixmier.
Thus
$$ \Afr\cong\smdp\Cpx{\alpha+\beta_2-1}{p\wedge q_2}\oplus
\Bigl((C([a,b])*A_0)\otimes M_2(\Cpx)\Bigr)
\oplus\smdp\Cpx{\alpha+\beta_1-1}{p\wedge q_1}. $$
Letting $r_0=1-p\wedge q_1-p\wedge q_2$ we then have that
$\Afr_0\eqdef r_0\Afr=(C([a,b])*A_0)\otimes M_2(\Cpx)$ is simple.
If $\phi\restrict_{A_0}$ is a trace
then $(1-p)\Afr(1-p)$ has stable rank~1 by
Proposition~\FPHaarUnitarySRone.
Thus
also $\Afr$ has stable rank~1.
Finally, clearly $r_0\Afr_1$ is in the centralizer of $\phi$.
Hence the centralizer of
$\phi\restrict_{\Afr_0}$ has a diffuse abelian subalgebra which contains
$r_0p$ and another which contains $r_0q_1$.

\proclaim{Case II(\AzeroCCC)} $\alpha=\beta_1>\frac12$. \endproclaim
Then
$$ \Afr_1=\{f:[0,b]\to M_2(\Cpx)\mid f\text{ continuous and }f(0)
\text{ diagonal }\}\oplus\smdp\Cpx{\alpha+\beta-1}{p\wedge q_1},
\tag{\CaseIIAone} $$
with $p=\left(\smallmatrix1&0\\0&0\endsmallmatrix\right)\oplus1$ and
$q_1=\left(\smallmatrix t&\sqrt{t(1-t)}\\
\sqrt{t(1-t)}&1-t\endsmallmatrix\right)\oplus1$.
Moreover, $p\wedge q_1$ is minimal and central in $\Afr$ and, by
Proposition~\FPHaarUnitaryDixmier,
$(1-p)\Afr(1-p)\cong C([0,b])*A_0$ is simple.
Consider the central projection, $r_0=1-p\wedge q_1$ and let $\Afr_0=r_0\Afr$.
Because $r_0\Afr_1$ is in the centralizer of $\phi$, the centralizer of
$\phi\restrict_{\Afr_0}$ has a diffuse abelian subalgebra which contains
$r_0p$ and and another which contains $r_0q_1$.
Let $\pi^{(1)}_{p\wedge q_2}:r_0\Afr_1\to\Cpx$ be the $*$--homomorphism
defined, in the
notation of~(\CaseIIAone), by
$$ \pi^{(1)}_{p\wedge q_2}(f)=\text{ the (1,1)--entry of }f(0), $$
so that $\pi^{(1)}_{p\wedge q_2}(r_0p)=1=\pi^{(1)}_{p\wedge q_2}(q_2)$.
Clearly $r_0$ is also central in $\Afr$ and the linear span of
$$ r_0p\Afr_1p+(1-p)\Afr(1-p)+(1-p)\Afr(1-p)\Afr_1p+p\Afr_1(1-p)\Afr(1-p)
+p\Afr_1(1-p)\Afr(1-p)\Afr_1p $$
is dense in $\Afr_0$.
Thus $\pi^{(1)}_{p\wedge q_2}$ extends to a
$*$--homomorphism $\pi_{p\wedge q_2}:\Afr_0\to\Cpx$ such that
$$ \aligned
\ker\pi^{(1)}_{p\wedge q_2}+
&(1-p)\Afr(1-p)+(1-p)\Afr(1-p)\Afr_1p+\\
+&p\Afr_1(1-p)\Afr(1-p)+p\Afr_1(1-p)\Afr(1-p)\Afr_1p
\endaligned \tag{\CaseIIdenseker} $$
spans a dense subset of $\ker\pi_{p\wedge q_2}$.

We now show that $\ker\pi_{p\wedge q_2}$ is simple.
Since $(1-p)\Afr(1-p)$ is simple, by Proposition~\FullSimple{} it will suffice
to show that $(1-p)\Afr(1-p)$ is full in $\ker\pi_{p\wedge q_2}$.
But clearly $1-p$ is full in $\ker\pi^{(1)}_{p\wedge q_2}$ and
$p\Afr_1(1-p)\subseteq\ker\pi^{(1)}_{p\wedge q_2}$.
Hence by the denseness of
the span of~(\CaseIIdenseker) in $\ker\pi_{p\wedge q_2}$, there is no proper
ideal of $\ker\pi_{p\wedge q_2}$ containing $(1-p)\Afr(1-p)$.

Suppose $\phi\restrict_{A_0}$ is a trace.
Then $(1-p)\Afr(1-p)$ has
stable rank~1 by Proposition~\FPHaarUnitarySRone.
Since $1-p$ is full in $\ker\pi_{p\wedge q_2}$,
also
$\ker\pi_{p\wedge q_2}$ has
stable rank~1 by Proposition~\AheredBsrone(i).
Then $\Afr$ has stable rank~1 by Proposition~\ShortExactSeq.

Finally, we show that $q_2$ is full in $\Afr_0$.
Suppose $\Ic$ is an ideal of $\Afr_0$ containing $q_2$.
Looking at the ideal of $\Afr_1$ generated by $q_2$, we see that $\Ic$ contains
a nonzero element of $\ker\pi_{p\wedge q_2}$, hence by simplicity
contains all of $\ker\pi_{p\wedge q_2}$.
But $q_2\not\in\ker\pi_{p\wedge q_2}$ and $\Afr_0/\ker\pi_{p\wedge q_2}$ is
one-dimensional, hence $\Afr_0\subseteq\Ic$.

\proclaim{Case III(\AzeroCCC)} $\beta_1>\alpha>\beta_2$. \endproclaim
Then
$$ \Afr_1=\smdp\Cpx{\beta_1-\alpha}{q_1\wedge(1-p)}\oplus
\Bigl(C([a,b])\otimes M_2(\Cpx)\Bigr)\oplus
\smdp\Cpx{\beta_1+\alpha-1}{q_1\wedge p} $$
and
$$ (1-p)\Afr(1-p)
\cong\Bigl(\smdp\Cpx{\frac{\beta_1-\alpha}{1-\alpha}}{q_1\wedge(1-p)}\oplus
C([a,b])\Bigr)*A_0 $$
is by Proposition~\FPHaarUnitaryDixmier{} simple.
Let $r_0=1-q_1\wedge p$ and $\Afr_0=r_0\Afr$.
Clearly $1-p$ is full in $r_0\Afr_1$ and thus is full in $\Afr_0$.
Hence $\Afr_0$ is simple.
If $D$ is a diffuse abelian subalgebra of the centralizer of
$\phi\restrict_{A_0}$ then $D+(p-q_1\wedge p)\Afr_1(p-q_1\wedge p)$ is a
diffuse abelian subalgebra of the centralizer of $\phi\restrict_{\Afr_0}$ and
contains $r_0p$.
By Proposition~\ExistsDAS{} and considering the C$^*$--subalgebra of
$(1-p)\Afr(1-p)$ generated by $\{q_1\wedge(1-p)\}\cup A_0$, we see that there
is a diffuse abelian subalgebra, $D$, of the centralizer of
$\phi\restrict_{(q_1\wedge(1-p))\Afr(q_1\wedge(1-p))}$.
Then
$$ D+q_2\Afr_1q_2
+(q_1-q_1\wedge(1-p)-q_1\wedge p)\Afr_1(q_1-q_1\wedge(1-p)-q_1\wedge p) $$
is a diffuse abelian subalgebra of the centralizer of $\phi\restrict_{\Afr_0}$
and contains $r_0q_1$.

If $\phi\restrict_{A_0}$ is a trace then
by Proposition~\FPHaarUnitarySRone{} $(1-p)\Afr(1-p)$ has stable rank~1.
So
by~Proposition~\AheredBsrone(i), $\Afr_0$ has stable rank~1.
We know
$$ \Afr=\Afr_0\oplus\smdp\Cpx{\alpha+\beta_1-1}{p\wedge q_1}, $$
so $\Afr$ has stable rank~1.

\proclaim{Case IV(\AzeroCCC)} $\frac12>\alpha=\beta_2$. \endproclaim
Then
$$ \Afr_1=\smdp\Cpx{\beta_1-\alpha}{(1-p)\wedge q_1}\oplus
\{f:[a,1]\to M_2(\Cpx)\mid f\text{ continuous and }f(1)\text{ diagonal }\},
\tag{\CaseIVAone} $$
with $0<a<1$,  $p=0\oplus\left(\smallmatrix1&0\\0&0\endsmallmatrix\right)$,
$q_1=1\oplus\left(\smallmatrix t&\sqrt{t(1-t)}\\
\sqrt{t(1-t)}&1-t\endsmallmatrix\right)$.
Thus
$$ (1-p)\Afr(1-p)\cong
\Bigl(\smdp\Cpx{\frac{\beta_1-\alpha}{1-\alpha}}{q_1\wedge(1-p)}\oplus
C([a,1])\Bigr)*A_0 $$
is by Proposition~\FPHaarUnitaryDixmier{} simple.
Moreover, if $D$ is a diffuse abelian subalgebra of the centralizer of
$\phi\restrict_{A_0}$ then $D+p\Afr_1p$ is a diffuse abelian subalgebra of the
centralizer of $\phi$ and
contains $p$.
By Proposition~\ExistsDAS{} and considering the C$^*$--subalgebra of
$(1-p)\Afr(1-p)$ generated by $\{q_1\wedge(1-p)\}\cup A_0$, we see that there
is a diffuse abelian subalgebra, $D$, of the centralizer of
$\phi\restrict_{(q_1\wedge(1-p))\Afr(q_1\wedge(1-p))}$.
Then
$$ D+q_2\Afr_1q_2+(q_1-q_1\wedge(1-p))\Afr_1(q_1-q_1\wedge(1-p)) $$
is a diffuse abelian subalgebra of the centralizer of $\phi$
and contains $q_1$.

Let $\pi^{(1)}_{p\wedge q_1}:\Afr_1\to\Cpx$ be the $*$--homomorphism defined,
in the notation of~(\CaseIVAone) by
$$ \pi^{(1)}_{p\wedge q_1}(\lambda\oplus f)=\text{ the (1,1)--entry of }f(1),
$$
so that $\pi^{(1)}_{p\wedge q_1}(p)=1=\pi^{(1)}_{p\wedge q_1}(q_1)$.
Then the linear span of
$$ p\Afr_1p+(1-p)\Afr(1-p)+(1-p)\Afr(1-p)\Afr_1p+p\Afr_1(1-p)\Afr(1-p)+
p\Afr_1(1-p)\Afr(1-p)\Afr_1p $$
is clearly dense in $\Afr$ so $\pi^{(1)}_{p\wedge q_1}$ extends to a
$*$--homomorphism $\pi_{p\wedge q_1}:\Afr\to\Cpx$ such that
$$ \aligned
\ker\pi^{(1)}_{p\wedge q_1}+
&(1-p)\Afr(1-p)+(1-p)\Afr(1-p)\Afr_1p+\\
+&p\Afr_1(1-p)\Afr(1-p)+p\Afr_1(1-p)\Afr(1-p)\Afr_1p
\endaligned \tag{\CaseIVdenseker} $$
spans a dense subset of $\ker\pi_{p\wedge q_1}$.

Like in Case~II, since $1-p$ is full in $\ker_{p\wedge q_1}^{(1)}$ and since
$(1-p)\Afr(1-p)$ is simple, it follows that $\ker\pi_{p\wedge q_1}$ is simple.

If $\phi\restrict_{A_0}$ is a trace then by Proposition~\FPHaarUnitarySRone,
$(1-p)\Afr(1-p)$ has
stable rank~1.
Since $1-p$ is full in $\ker\pi_{p\wedge q_1}$,
also $\ker\pi_{p\wedge q_1}$ has
stable rank~1 by
Proposition~\AheredBsrone(i).
Thus by Proposition~\ShortExactSeq, $\Afr$ has stable rank~1.

Finally, we show that $q_1$ is full in $\Afr$.
Suppose $\Ic$ is an ideal of $\Afr$ containing $q_1$.
Looking at the ideal of $\Afr_1$ generated by $q_1$, we see that $\Ic$ contains
a nonzero element of $\ker\pi_{p\wedge q_1}$, hence by simplicity
contains all of $\ker\pi_{p\wedge q_1}$.
But $q_1\not\in\ker\pi_{p\wedge q_1}$ and $\pi_{p\wedge q_1}$ is
one-dimensional, hence $\Afr\subseteq\Ic$.

\proclaim{Case V(\AzeroCCC)} $\beta_2>\alpha$. \endproclaim
Then
$$ \Afr_1=\smdp\Cpx{\beta_1-\alpha}{(1-p)\wedge q_1}\oplus
\Bigl(C([a,b])\otimes M_2(\Cpx)\Bigr)\oplus
\smdp\Cpx{\beta_2-\alpha}{(1-p)\wedge q_2} $$
and
$$ (1-p)\Afr(1-p)\cong\Bigl(
\smdp\Cpx{\frac{\beta_1-\alpha}{1-\alpha}}{(1-p)\wedge q_1}\oplus
C([a,b])\oplus
\smdp\Cpx{\frac{\beta_2-\alpha}{1-\alpha}}{(1-p)\wedge q_2}\Bigr)*A_0 $$
is by Proposition~\FPHaarUnitaryDixmier{} simple.
Since $1-p$ is full in $\Afr_1$, is follows that $\Afr$ is simple.

Moreover, if $D$ is a diffuse abelian subalgebra of the centralizer of
$\phi\restrict_{A_0}$ then $D+p\Afr_1p$ is a diffuse abelian subalgebra of
the centralizer of $\phi$ and
contains $p$.
By Proposition~\ExistsDAS{} and considering the C$^*$--subalgebra of
$(1-p)\Afr(1-p)$ generated by $\{q_1\wedge(1-p),q_2\wedge(1-p)\}\cup A_0$, we
see that there are
a diffuse abelian subalgebras, $D_1$ and $D_2$, of the centralizers of
$\phi\restrict_{(q_1\wedge(1-p))\Afr(q_1\wedge(1-p))}$ and, respectively,
$\phi\restrict_{(q_2\wedge(1-p))\Afr(q_2\wedge(1-p))}$.
Then
$$ D_1+(q_1-q_1\wedge(1-p))\Afr_1(q_1-q_1\wedge(1-p))+D_2
+(q_2-q_2\wedge(1-p))\Afr_1(q_2-q_2\wedge(1-p)) $$
is a diffuse abelian subalgebra of the centralizer of $\phi$
and contains $q_1$.

If $\phi\restrict_{A_0}$ is a trace then by Proposition~\FPHaarUnitarySRone{}
$(1-p)\Afr(1-p)$ has stable rank~1.
Hence by Proposition~\AheredBsrone(i) also $\Afr$ has stable rank~1.

\proclaim{Case VI(\AzeroCCC)} $\beta_1=\alpha=\frac12$. \endproclaim
Then
$$ \Afr_1=\{f:[0,1]\to M_2(\Cpx)\mid f\text{ continuous, }f(0),\,f(1)
\text{ diagonal }\}, \tag{\CaseVIAone} $$
with $p=\left(\smallmatrix1&0\\0&0\endsmallmatrix\right)$,
$q_1=\left(\smallmatrix t&\sqrt{t(1-t)}\\
\sqrt{t(1-t)}&1-t\endsmallmatrix\right)$.
Thus
$$ (1-p)\Afr(1-p)\cong C([0,1])*A_0 $$
is by Proposition~\FPHaarUnitaryDixmier{} simple.
Moreover, clearly $\Afr_1$ is in the centralizer of $\phi$ and has diffuse
abelian subalgebras containing $p$ and, respectively, $q_1$.
For $i\in\{1,2\}$
let $\pi^{(1)}_{p\wedge q_i}:\Afr_1\to\Cpx$ be the $*$--homomorphism defined,
in the notation of~(\CaseVIAone), by
$$ \pi^{(1)}_{p\wedge q_i}(f)=
\cases\text{ the (1,1)--entry of }f(1)&\text{if }i=1\\
\text{ the (1,1)--entry of }f(0)&\text{if }i=2, \endcases $$
so that $\pi^{(1)}_{p\wedge q_i}(p)=1=\pi^{(1)}_{p\wedge q_i}(q_i)$.
Then the linear span of
$$ p\Afr_1p+(1-p)\Afr(1-p)+(1-p)\Afr(1-p)\Afr_1p+p\Afr_1(1-p)\Afr(1-p)+
p\Afr_1(1-p)\Afr(1-p)\Afr_1p $$
is clearly dense in $\Afr$ so $\pi^{(1)}_{p\wedge q_i}$ extends to a
$*$--homomorphism $\pi_{p\wedge q_i}:\Afr\to\Cpx$ such that
$$ \aligned
\ker\pi^{(1)}_{p\wedge q_1}\cap\ker\pi^{(1)}_{p\wedge q_2}+
&(1-p)\Afr(1-p)+(1-p)\Afr(1-p)\Afr_1p+\\
+&p\Afr_1(1-p)\Afr(1-p)+p\Afr_1(1-p)\Afr(1-p)\Afr_1p
\endaligned \tag{\CaseVIdenseker} $$
spans a dense subset of $\ker\pi_{p\wedge q_1}\cap\ker\pi_{p\wedge q_2}$.

Like in Case~II, since $1-p$ is full in
$\ker_{p\wedge q_1}^{(1)}\cap\ker_{p\wedge q_2}^{(1)}$
and since $(1-p)\Afr(1-p)$ is simple,
it follows that $\ker\pi_{p\wedge q_1}\cap\ker\pi_{p\wedge q_2}$ is simple.

If $\phi\restrict_{A_0}$ is a trace then by Proposition~\FPHaarUnitarySRone{}
$(1-p)\Afr(1-p)$ has stable rank~1.
Since $1-p$ is full in $\ker\pi_{p\wedge q_1}\cap\ker\pi_{p\wedge q_2}$,
by Proposition~\AheredBsrone(i) also
$\ker\pi_{p\wedge q_1}\cap\ker\pi_{p\wedge q_2}$ has stable rank~1.
Since $\Afr/(\ker\pi_{p\wedge q_1}\cap\ker\pi_{p\wedge q_2})$ is
two--dimensional, it follows from Proposition~\ShortExactSeq{} that $\Afr$ has
stable rank~1.

We show that $q_1$ is full in $\ker\pi_{p\wedge q_2}$.
Suppose $\Ic$ is an ideal of $\ker\pi_{p\wedge q_2}$ containing $q_1$.
Multiplying by elements of $\Afr_1$, we see that $\Ic$ contains a nonzero
element of $\ker\pi_{p\wedge q_1}$, hence by simplicity
contains all of $\ker\pi_{p\wedge q_1}\cap\ker\pi_{p\wedge q_2}$.
But $q_1\not\in\ker\pi_{p\wedge q_1}$ and $\pi_{p\wedge q_1}$ is
one-dimensional, hence $\ker\pi_{p\wedge q_2}\subseteq\Ic$.
The proof that
$q_2$ is full in $\ker\pi_{p\wedge q_1}$
is the same.

We now examine the question of existence and uniqueness of tracial states on
the algebras delineated in the statement of the lemma.
In all the cases above, it follows from
Proposition~\FPHaarUnitaryDixmier{} that $(1-p)\Afr(1-p)$ has tracial states if
and only if $\phi\restrict_{A_0}$ is a trace, and then the free product state
gives the unique tracial state on $(1-p)\Afr(1-p)$.
Moreover, the element $1-p$ is full in the simple algebras under consideration,
i.e.
\roster
\item"" $1-p$ is full in $\Afr_0$ in Cases~I,~III and~V,
\item"" $1-p$ is full in $\ker\pi_{p\wedge q_2}$ in Case~II,
\item"" $1-p$ is full in $\ker\pi_{p\wedge q_1}$ in Case~IV,
\item"" $1-p$ is full in $\ker\pi_{p\wedge q_1}\cap\ker\pi_{p\wedge q_2}$ in
Case~VI.
\endroster
It then follows from
Proposition~\AheredBsrone(ii) that in each of the Cases~I--VI, the
corresponding algebra has tracial states if
and only if $\phi\restrict_{A_0}$ is a trace, and then the restriction of
$\phi(r_0)^{-1}\phi$ to this algebra is
its unique tracial state.
(In the non--unital Cases~II,~IV and~VI,
one easily sees that the above normalization gives a state by looking at the
subalgebra $r_0\Afr_1$.)
\QED

\proclaim{Lemma \AzeroCCCC}
Let
$$ (\Afr,\phi)=(A_0\oplus\smdp\Cpx{\alpha_1}{p_1}\oplus
\smdp\Cpx{\alpha_2}{p_2})*
(\smdp\Cpx{\beta_1}{q_1}\oplus\smdp\Cpx{\beta_2}{q_2}), $$
where the centralizer of $\phi\restrict_{A_0}$ has a diffuse abelian
subalgebra.
Let
$$ \aligned
L_+&=\{(i,j)\mid\alpha_i+\beta_j>1\}\\
L_0&=\{(i,j)\mid\alpha_i+\beta_j=1\}.
\endaligned \tag{\LplusLzero} $$
Then
$$ \Afr=\smp{\Afr_0}{r_0}\oplus
\bigoplus_{(i,j)\in L_+}\smdp\Cpx{\alpha_i+\beta_j-1}{p_i\wedge q_j}, $$
where the centralizer of $\phi\restrict_{\Afr_0}$ has a diffuse abelian
subalgebra which contains $r_0p_1$ and a diffuse abelian subalgebra
which contains $r_0q_1$.

If $\phi\restrict_{A_0}$ is a trace then the stable rank of $\Afr$ is~1.

If $L_0$ is empty then $\Afr_0$ is simple.
If, in addition, $\phi\restrict_{A_0}$ is a trace then
$\phi(r_0)^{-1}\phi\restrict_{\Afr_0}$
is the unique tracial state
on $\Afr_0$ and if $\phi\restrict_{A_0}$ is not a trace then $\Afr_0$ has no
tracial states.

If $L_0$ is not empty, then for every $(i,j)\in L_0$ there is a
$*$--homomorphism $\pi_{(i,j)}:\Afr_0\to\Cpx$ such that
$\pi_{(i,j)}(r_0p_i)=1=\pi_{(i,j)}(r_0q_j)$.
Then
\roster
\item"(i)"
$$ \Afr_{00}\eqdef\dsize\bigcap_{(i,j)\in L_0}\ker\pi_{(i,j)} $$
is simple.
If $\phi\restrict_{A_0}$ is a trace then
$\phi(r_0)^{-1}\phi\restrict_{\Afr_{00}}$
is the unique tracial state
on $\Afr_{00}$ and if $\phi\restrict_{A_0}$ is not a trace then $\Afr_{00}$ has no
tracial states.
\item"(ii)" For each $i\in\{1,2\}$, $r_0p_i$ is full in
$\dsize\Afr_0\cap\bigcap\Sb(i',j)\in L_0\\i'\ne i\endSb \ker\pi_{(i',j)}$.
\vskip2ex
\item"(iii)" For each $j\in\{1,2\}$, $r_0q_j$ is full in
$\dsize\Afr_0\cap\bigcap\Sb(i,j')\in L_0\\j'\ne j\endSb \ker\pi_{(i,j')}$.
\endroster
\endproclaim
\demo{Proof}
We will assume that $\alpha_1\ge\alpha_2$ and $\beta_1\ge\beta_2$.
To prove the lemma in its full generality, we will now be careful to find a
diffuse abelian subalgebra of the centralizer of $\phi\restrict_{\Afr_0}$
containing $\{r_0p_1,r_0p_2\}$, not only $r_0p_1$.
Let $\Afr_1$ be the C$^*$--subalgebra of $\Afr$ generated by
$A_0+\Cpx(p_1+p_2)$ together with $\{q_1,q_2\}$, i.e\.
$$ (\Afr_1,\phi\restrict_{\Afr_1})=
(A_0\oplus\smdp\Cpx{\alpha_1+\alpha_2}{p_1+p_2})*
(\smdp\Cpx{\beta_1}{q_1}\oplus\smdp\Cpx{\beta_2}{q_2}).
\tag{\AzeroCCCCAfrone} $$
We find $\Afr_1$ using Lemma~\AzeroCCC.
Then, by Proposition~\FreeProdDirectSum
$$ (p_1+p_2)\Afr(p_1+p_2)\cong(p_1+p_2)\Afr_1(p_1+p_2)*
\Bigl(\smdp\Cpx{\frac{\alpha_1}{\alpha_1+\alpha_2}}{p_1}\oplus
\smdp\Cpx{\frac{\alpha_2}{\alpha_1+\alpha_2}}{p_2}\Bigr) \tag{\AzeroCCCCAfr} $$
We consider three cases.

\proclaim{Case I(\AzeroCCCC)} $\alpha_1+\alpha_2+\beta_1\le1$. \endproclaim
Then by Lemma~\AzeroCCC, the centralizer of $\phi\restrict_{\Afr_1}$ has a
diffuse abelian subalgebra, $D$, which contains $p_1+p_2$ and
$p_1+p_2$ is full in $\Afr_1$.
Hence by Proposition~\FPHaarUnitaryDixmier{} $(p_1+p_2)\Afr(p_1+p_2)$ is
simple.
Also, $p_1+p_2$ is full in $\Afr$, hence $\Afr$ is simple.
If $\phi_{A_0}$ is a trace then by Proposition~\FPHaarUnitarySRone{}
$(p_1+p_2)\Afr(p_1+p_2)$ has stable rank~1.
Hence by Proposition~\AheredBsrone(i) so does $\Afr$.
The application of Lemma~\AzeroCCC{} to~(\AzeroCCCCAfrone) yields a diffuse
abelian subalgebra of the centralizer of $\phi$ which contains $\{q_1,q_2\}$.
Applying Corollary~\ExistsDAScorr{} to~(\AzeroCCCCAfr) shows that the
centralizer
of $\phi\restrict_{(p_1+p_2)\Afr(p_1+p_2)}$ has a diffuse abelian subalgebra,
$D'$, containing $\{p_1,p_2\}$.
Then $(1-p_1-p_2)D(1-p_1-p_2)+D'$ is a diffuse abelian subalgebra of the
centralizer of $\phi$ containing $\{p_1,p_2\}$.

\proclaim{Case II(\AzeroCCCC)} $\alpha_1+\alpha_2+\beta_1>1$ and
$\alpha_1+\alpha_2+\beta_2\le1$. \endproclaim
Note that this implies
$$ \aligned
L_+&=\{(i,1)\mid\alpha_i+\beta_1>1\}\\
L_0&=\{(i,1)\mid\alpha_i+\beta_1=1\}.
\endaligned $$
Applying Lemma~\AzeroCCC{} to~(\AzeroCCCCAfrone) shows that
$$ \Afr_1=\smp{\Afr_{1,0}}{r_{1,0}}\oplus
\smdp\Cpx{\alpha_1+\alpha_2+\beta_1-1}{(p_1+p_2)\wedge q_1}, $$
where $r_{1,0}=1-(p_1+p_2)\wedge q_1$, where
the centralizer of $\phi\restrict_{\Afr_{1,0}}$ has a diffuse abelian
subalgebra containing $r_{1,0}(p_1+p_2)$ and
where each of $r_{1,0}(p_1+p_2)$ and $q_2$ is full in $\Afr_{1,0}$.
Then
$$ (p_1+p_2)\Afr_1(p_1+p_2)=(p_1+p_2)\Afr_{1,0}(p_1+p_2)\oplus
\smdp\Cpx{\frac{\alpha_1+\alpha_2+\beta_1-1}{\alpha_1+\alpha_2}}
{(p_1+p_1)\wedge q_1}. $$
So from~(\AzeroCCCCAfr) and Lemma~\AzeroCCC,
$$ (p_1+p_2)\Afr(p_1+p_2)
=\smp{\Afr_{2,0}}{r_{2,0}}\oplus\bigoplus_{(i,1)\in L_+}
\smdp\Cpx{\frac{\alpha_i+\beta_1-1}{\alpha_1+\alpha_2}}{p_i\wedge q_1}, $$
where $r_{2,0}=p_1+p_2-\sum_{(i,1)\in L_+}p_i\wedge q_1$.
Since $\Afr$ is generated by $(p_1+p_2)\Afr(p_1+p_2)\cup\Afr_1$, we have
$$ \Afr=\smp{\Afr_0}{r_0}\oplus\bigoplus_{(i,1)\in L_+}
\smdp\Cpx{\alpha_i+\beta_1-1}{p_i\wedge q_1} $$
where the linear span of
$$ \aligned
\Afr_{2,0}&+\Afr_{2,0}\Afr_{1,0}(1-p_1-p_2)+(1-p_1-p_2)\Afr_{1,0}\Afr_{2,0}+\\
&+(1-p_1-p_2)\Afr_{1,0}\Afr_{2,0}\Afr_{1,0}(1-p_1-p_2)+
(1-p_1-p_2)\Afr_{1,0}(1-p_1-p_2) 
\endaligned \tag{\CaseIIAzeroCCCCdense} $$
is dense in $\Afr_0$.
Thus $r_0=r_{2,0}+(1-p_1-p_2)$.
Now the centralizer of
$\phi\restrict_{\Afr_{2,0}}$ has a diffuse abelian subalgebra, $D$, which
contains $\{r_{2,0}p_1,r_{2,0}p_2\}$.
Letting $D'$ be a diffuse abelian subalgebra of the centralizer of
$\phi\restrict_{A_0}$ it follows that $D+D'$ is a diffuse abelian subalgebra of
the centralizer of $\phi\restrict_{\Afr_0}$ and contains $\{r_0p_1,r_0p_2\}$.

Let now $D$ be a diffuse abelian subalgebra of the centralizer of
$\phi\restrict_{\Afr_{2,0}}$ containing $r_{2,0}((p_1+p_2)\wedge q_1)$
and $D'$ be
a diffuse abelian subalgebra of the centralizer of $\phi\restrict_{\Afr_{1,0}}$
containing $\{r_{1,0}q_1,q_2\}$.
Then
$r_{2,0}((p_1+p_2)\wedge q_1)D+D'$
is a diffuse abelian subalgebra of the centralizer of $\phi\restrict_{\Afr_0}$
and contains $\{q_2,r_0q_1\}$.

Since $r_{1,0}(p_1+p_2)\in\Afr_{2,0}$ and is full in $\Afr_{1,0}$, it follows
that $\Afr_{2,0}$ is full in $\Afr_0$.
If $L_0$ is empty then $\Afr_{2,0}$ is simple, hence (by
Proposition~\FullSimple) $\Afr_0$ is also simple.

Otherwise, if $L_0$ is nonempty, for every $(i,1)\in L_0$ there is a
$*$--homomorphism $\pi_{(i,1)}^{(2)}:\Afr_{2,0}\to\Cpx$ such that
$\pi_{(i,1)}^{(2)}(r_{2,0}p_i)=1
=\pi_{(i,1)}^{(2)}(r_{2,0}((p_1+p_2)\wedge q_1))$.
Using the denseness of the span of~(\CaseIIAzeroCCCCdense) in $\Afr_0$, we see
that
$\pi_{(i,1)}^{(2)}$ extends to a $*$--homomorphism $\pi_{(i,1)}:\Afr_0\to\Cpx$
such that
$\pi_{(i,1)}(r_0p_i)=1=\pi_{(i,1)}(r_0q_1)$ and the linear span of
$$ \align
\ker\pi_{(i,1)}^{(2)}&+\Afr_{2,0}\Afr_{1,0}(1-p_1-p_2)+
(1-p_1-p_2)\Afr_{1,0}\Afr_{2,0}+\\
&+(1-p_1-p_2)\Afr_{1,0}\Afr_{2,0}\Afr_{1,0}(1-p_1-p_2)+
(1-p_1-p_2)\Afr_{1,0}(1-p_1-p_2) 
\endalign $$
is dense $\ker\pi_{(i,1)}$.

Let
$$ \Afr_{2,00}=\Afr_{2,0}\cap\bigcap_{(i,1)\in L_0}\ker\pi_{(i,1)}^{(2)}. $$
From the application of Lemma~\AzeroCCC, $\Afr_{2,00}$ is simple.
Since $\Afr_{2,00}$ contains $r_{1,0}(p_1+p_2)$, which is full in $\Afr_{1,0}$,
it follows that $\Afr_{2,00}$ is full in $\Afr_{00}$.
Then (by Proposition~\FullSimple), $\Afr_{00}$ is simple.

Let $i\in\{1,2\}$.
We now show that $r_0p_i$ is full in
$$ \Afr_0\cap\dsize\bigcap\Sb (i',1)\in L_0\\i'\ne i\endSb\ker\pi_{(i',1)}. 
\tag{\CaseIIAzeroCCCCfullin} $$
Suppose $\Ic$ is an ideal of the algebra in~(\CaseIIAzeroCCCCfullin) containing
$r_0p_i$.
Since $r_{2,0}\le r_0$ and (by Lemma~\AzeroCCC) $r_{2,0}p_i$ is full in
$$ \Afr_{2,0}\cap\dsize\bigcap\Sb (i',1)\in L_0\\i'\ne i\endSb
\ker\pi_{(i',1)}^{(2)}, \tag{\CaseIIAzeroCCCCfullintwo} $$
$\Ic$ must contain the algebra in~(\CaseIIAzeroCCCCfullintwo).
Hence, arguing as above, $\Afr_{1,0}\subseteq\Ic$.
Thus
$$ \align
\Afr_{2,0}\cap\dsize\bigcap\Sb (i',1)\in L_0\\i'\ne i\endSb
\ker&\pi_{(i',1)}^{(2)}+
\Afr_{2,0}\Afr_{1,0}(1-p_1-p_2)+\\
&+(1-p_1-p_2)\Afr_{1,0}\Afr_{2,0}
+(1-p_1-p_2)\Afr_{1,0}\Afr_{2,0}\Afr_{1,0}(1-p_1-p_2)\subseteq\Ic,
\endalign $$
proving that the algebra of~(\CaseIIAzeroCCCCfullin) is contained in $\Ic$.

Similarly, since $r_{2,0}((p_1+p_2)\wedge q_1)$ is full in $\Afr_{2,0}$, it
follows that $r_0q_1$
is full in $\Afr_0$.

If $\phi\restrict_{A_0}$ is a trace then from Lemma~\AzeroCCC{} we have that
$\Afr_2$ and indeed $\Afr_{2,0}$ has stable rank~1.
The fullness of $\Afr_{2,0}$ in $\Afr_0$ implies (via
Proposition~\AheredBsrone(i)) that $\Afr_0$ has stable rank~1, hence $\Afr$ has
stable rank~1.

Finally, concerning existence and uniqueness of traces,
from Lemma~\AzeroCCC{} we have that $\Afr_{2,00}$ has a tracial state if and
only if $\phi\restrict_{A_0}$ is a trace, and then
$\phi(r_{2,0})^{-1}\phi\restrict_{\Afr_{2,00}}$ is its unique tracial state.
The same statement for $\Afr_{00}$ then follows from fullness of
$\Afr_{2,00}$ in $\Afr_{00}$ and Proposition~\AheredBsrone(ii).
(One can easily check the normalization.)

\proclaim{Case III(\AzeroCCCC)} $\alpha_1+\alpha_2+\beta_2>1$. \endproclaim
Since $\beta_2\le\frac12$ we must have $\frac12<\alpha_1+\alpha_2$.
Let $n\in\Naturals$ be least such that $\alpha_1+\alpha_2\le\frac n{n+1}$.
Thus $n\ge2$.
We will proceed by induction on $n$, proving the case $n=2$ and the inductive
step simultaneously.
Applying Lemma~\AzeroCCC{} to~(\AzeroCCCCAfrone) we have
$$ \Afr_1=\smp{\Afr_{1,0}}{r_{1,0}}\oplus
\smdp\Cpx{\alpha_1+\alpha_2+\beta_1-1}{(p_1+p_2)\wedge q_1}\oplus
\smdp\Cpx{\alpha_1+\alpha_2+\beta_2-1}{(p_1+p_2)\wedge q_2}, $$
where $r_{1,0}=1-(p_1+p_2)\wedge q_1-(p_1+p_2)\wedge q_2$,
where the centralizer of $\phi\restrict_{\Afr_{1,0}}$ has a diffuse abelian
subalgebra containing $r_{1,0}(p_1+p_2)$ and where $\Afr_{1,0}$ is simple.
Thus from~(\AzeroCCCCAfr),
$$ (p_1+p_2)\Afr(p_1+p_2)\cong
\Bigl((p_1+p_2)\Afr_{1,0}(p_1+p_2)\oplus
\smdp\Cpx{\frac{\alpha_1+\alpha_2+\beta_1-1}{\alpha_1+\alpha_2}}
{(p_1+p_2)\wedge q_1}\oplus
\smdp\Cpx{\frac{\alpha_1+\alpha_2+\beta_2-1}{\alpha_1+\alpha_2}}
{(p_1+p_2)\wedge q_2}\Bigr)*
\Bigl(\smdp\Cpx{\frac{\alpha_1}{\alpha_1+\alpha_2}}{p_1}\oplus
\smdp\Cpx{\frac{\alpha_2}{\alpha_1+\alpha_2}}{p_2}\Bigr). $$
Now since $\frac{n-1}n<\alpha_1+\alpha_2$, we have
$$ \frac{\alpha_1+\alpha_2+\beta_1-1}{\alpha_1+\alpha_2}+
\frac{\alpha_1+\alpha_2+\beta_2-1}{\alpha_1+\alpha_2}
=2-\frac1{\alpha_2+\alpha_2}<\frac{n-2}{n-1}. $$
The inductive hypothesis (or, when $n\in\{2,3\}$, the previously
considered Case~I or Case~II,) applies and we have, with $L_+$ as
in~(\LplusLzero),
$$ (p_1+p_2)\Afr(p_1+p_2)=\smp{\Afr_{2,0}}{r_{2,0}}\oplus
\bigoplus_{(i,j)\in L_+}
\smdp\Cpx{\frac{\alpha_i+\beta_j-1}{\alpha_1+\alpha_2}}{p_i\wedge q_j}, $$
where $r_{2,0}=p_1+p_2-\sum_{(i,j)\in L_+}p_i\wedge q_j$.
We obtain that
$$ \Afr=\smp{\Afr_0}{r_0}\oplus
\bigoplus_{(i,j)\in L_+}\smdp\Cpx{\alpha_i+\beta_j-1}{p_i\wedge q_j} $$
and that the span of~(\CaseIIAzeroCCCCdense) is dense in $\Afr_0$.
So $r_0=r_{2,0}+(1-p_1-p_2)$.
Letting $D$ be a diffuse abelian subalgebra of
the centralizer of $\phi\restrict_{\Afr_{2,0}}$ 
containing $\{r_{2,0}p_1,r_{2,0}p_2\}$
and $D'$ be a diffuse abelian subalgebra of the centralizer of
$\phi\restrict_{A_0}$, we see that $D+D'$ is a diffuse abelian subalgebra of
the centralizer of $\phi$ and contains $\{r_0p_1,r_0p_2\}$.

Let $D$ be a diffuse abelian subalgebra of the centralizer of
$\phi\restrict_{\Afr_{2,0}}$ that contains
$\{((p_1+p_2)\wedge q_1)r_{2,0},((p_1+p_2)\wedge q_2)r_{2,0}\}$ and let $D'$ be
a diffuse abelian subalgebra of the centralizer of $\phi\restrict_{\Afr_{1,0}}$
that contains $\{r_{1,0}q_1,r_{1,0}q_2\}$.
Then
$$ r_{2,0}((p_1+p_2)\wedge q_1)D+r_{2,0}((p_1+p_2)\wedge q_2)D+D' $$
is a diffuse abelian subalgebra of the centralizer of $\phi\restrict_{\Afr_0}$
and contains $\{r_0q_1,r_0q_2\}$.

Since $r_{1,0}(p_1+p_2)\in\Afr_{2,0}$ and is full in $\Afr_{1,0}$, it follows
that $\Afr_{2,0}$ is full in $\Afr_0$.
If $L_0$, as in~(\LplusLzero), is empty then $\Afr_{2,0}$ is simple, hence (by
Proposition~\FullSimple) $\Afr_0$ is also simple.

If $L_0$ is nonempty, then for every $(i,j)\in L_0$ there is a
$*$--homomorphism $\pi_{(i,j)}^{(2)}:\Afr_{2,0}\to\Cpx$ such that
$\pi_{(i,j)}^{(2)}(r_{2,0}p_i)=1
=\pi_{(i,j)}^{(2)}(r_{2,0}((p_1+p_2)\wedge q_j))$, and this extends to a
$*$--homomorphism $\pi_{(i,j)}:\Afr_0\to\Cpx$ such that
$\pi_{(i,j)}(r_0p_i)=1=\pi_{(i,j)}(r_0q_j)$ and the linear span of
$$ \align
\ker\pi_{(i,j)}^{(2)}&+\Afr_{2,0}\Afr_{1,0}(1-p_1-p_2)+
(1-p_1-p_2)\Afr_{1,0}\Afr_{2,0}+\\
&+(1-p_1-p_2)\Afr_{1,0}\Afr_{2,0}\Afr_{1,0}(1-p_1-p_2)+
(1-p_1-p_2)\Afr_{1,0}(1-p_1-p_2) 
\endalign $$
is dense $\ker\pi_{(i,j)}$.

Let
$$ \Afr_{2,00}=\Afr_{2,0}\cap\bigcap_{(i,j)\in L_0}\ker\pi_{(i,j)}^{(2)}. $$
From the application of Lemma~\AzeroCCC, $\Afr_{2,00}$ is simple.
Since $\Afr_{2,00}$ contains $r_{1,0}(p_1+p_2)$, which is full in $\Afr_{1,0}$,
it follows that $\Afr_{2,00}$ is full in $\Afr_{00}$.
Then (by Proposition~\FullSimple), $\Afr_{00}$ is simple.

Let $i\in\{1,2\}$.
We will now show that $r_0p_i$ is full in
$$ \Afr_0\cap\bigcap\Sb(i',j)\in L_0\\i'\ne i\endSb
\ker\pi_{(i',j)}. \tag{\CaseIIIAzeroCCCCfullini} $$
Suppose $\Ic$ is an ideal of the algebra in~(\CaseIIIAzeroCCCCfullini) which
contains $r_0p_i$.
Since $r_{2,0}\le r_0$ and since $r_{2,0}p_i$ is full in
$$ \Afr_{2,0}\cap\bigcap\Sb(i',j)\in L_0\\i'\ne i\endSb
\ker\pi_{(i',j)}^{(2)}, $$
this algebra must be contained in $\Ic$.
Then $r_{1,0}(p_1+p_2)\in\Ic$.
Since $\Afr_{1,0}$ is simple, it is then contained in $\Ic$.
Hence
$$ \align
\Afr_{2,0}\cap\dsize\bigcap\Sb (i',j)\in L_0\\i'\ne i\endSb
\ker&\pi_{(i',j)}^{(2)}+
\Afr_{2,0}\Afr_{1,0}(1-p_1-p_2)+\\
&+(1-p_1-p_2)\Afr_{1,0}\Afr_{2,0}
+(1-p_1-p_2)\Afr_{1,0}\Afr_{2,0}\Afr_{1,0}(1-p_1-p_2)\subseteq\Ic,
\endalign $$
proving that the algebra of~(\CaseIIIAzeroCCCCfullini) is contained in $\Ic$.

Let $j\in\{1,2\}$.
We now show that $r_0q_j$ is full in
$$ \Afr_0\cap\bigcap\Sb(i,j')\in L_0\\j'\ne j\endSb
\ker\pi_{(i,j')}. \tag{\CaseIIIAzeroCCCCfullinj} $$
Suppose $\Ic$ is an ideal of the algebra in~(\CaseIIIAzeroCCCCfullinj) which
contains $r_0q_j$.
Since $r_{2,0}\le r_0$ and since $r_{2,0}q_j$ is full in
$$ \Afr_{2,0}\cap\bigcap\Sb(i,j')\in L_0\\j'\ne j\endSb
\ker\pi_{(i,j')}^{(2)}, $$
this algebra must be contained in $\Ic$.
Then $r_{1,0}(p_1+p_2)\in\Ic$, which as before shows
that the algebra of~(\CaseIIIAzeroCCCCfullinj) is contained in $\Ic$.

The required results about the stable rank of $\Afr$ and the existence and
uniqueness of traces on $\Afr_{00}$ follow from the inductive hypothesis
because
(in the simple case) $\Afr_{2,0}$ is full in $\Afr_0$ or (more generally)
$\Afr_{2,00}$ is full in $\Afr_{00}$.
\QED

\proclaim{Lemma \CdCCC}
Let $n\in\Naturals$, $n\ge3$ and let
$$ (\Afr,\phi)=(\smdp\Cpx{\alpha_1}{p_1}\oplus\cdots\oplus
\smdp\Cpx{\alpha_n}{p_n})*
(\smdp\Cpx{\beta_1}{q_1}\oplus\smdp\Cpx{\beta_2}{q_2}). $$
Let
$$ \aligned
L_+&=\{(i,j)\mid\alpha_i+\beta_j>1\}\\
L_0&=\{(i,j)\mid\alpha_i+\beta_j=1\}.
\endaligned \tag{\LplusLzeron} $$
Then
$$ \Afr=\smp{\Afr_0}{r_0}\oplus
\bigoplus_{(i,j)\in L_+}\smdp\Cpx{\alpha_i+\beta_j-1}{p_i\wedge q_j}, $$
where the centralizer of $\phi\restrict_{\Afr_0}$ has a diffuse abelian
subalgebra which contains $r_0p_1$ and a diffuse abelian subalgebra
which contains $r_0q_1$.

The stable rank of $\Afr$ is~1.

If $L_0$ is empty then $\Afr_0$ is simple and
$\phi(r_0)^{-1}\phi\restrict_{\Afr_0}$
is the unique tracial state
on $\Afr_0$.

If $L_0$ is not empty, then for every $(i,j)\in L_0$ there is a
$*$--homomorphism $\pi_{(i,j)}:\Afr_0\to\Cpx$ such that
$\pi_{(i,j)}(r_0p_i)=1=\pi_{(i,j)}(r_0q_j)$.
Then
\roster
\item"(i)"
$$ \Afr_{00}\eqdef\dsize\bigcap_{(i,j)\in L_0}\ker\pi_{(i,j)} $$
is simple and
$\phi(r_0)^{-1}\phi\restrict_{\Afr_{00}}$
is the unique tracial state
on $\Afr_0$.
\item"(ii)" For each $i\in\{1,2,\ldots,n\}$, $r_0p_i$ is full in
$$ \Afr_0\cap\bigcap\Sb(i',j)\in L_0\\i'\ne i\endSb \ker\pi_{(i',j)}.
\tag{\CdCCCdensei} $$
\vskip2ex
\item"(iii)" For each $j\in\{1,2\}$, $r_0q_j$ is full in
$$ \Afr_0\cap\bigcap\Sb(i,j')\in L_0\\j'\ne j\endSb\ker\pi_{(i,j')}.
\tag{\CdCCCdensej} $$
\endroster
\endproclaim
\demo{Proof}
We proceed by induction on $n$, proving the initial step $n=3$ and the
inductive step simultaneously.
Let $\Afr_1$ be the C$^*$--subalgebra of $\Afr$ generated by
$(\Cpx(p_1+p_2)+\Cpx p_3+\cdots\Cpx p_n)\cup(\Cpx q_1+\Cpx q_2)$.
Thus
$$ (\Afr_1,\phi\restrict_{\Afr_1})\cong
(\smdp\Cpx{\alpha_1+\alpha_2}{p_1+p_2}\oplus\smdp\Cpx{\alpha_3}{p_3}
\oplus\cdots\oplus\smdp\Cpx{\alpha_n}{p_n})
*(\smdp\Cpx{\beta_1}{q_1}\oplus\smdp\Cpx{\beta_2}{q_2}). $$
By the inductive hypothesis when $n>3$ or by Proposition~\TwoProj{} when $n=3$,
letting
$$ \alignedat 2
L_+&=\{(i,j)\mid\alpha_i+\beta_j>1\}
&L_0&=\{(i,j)\mid\alpha_i+\beta_j=1\} \\
L_+^{(1)}&=\{(i,j)\mid i\ge3,\,\alpha_i+\beta_j>1\}
\qquad&L_0^{(1)}&=\{(i,j)\mid i\ge3,\,\alpha_i+\beta_j=1\} \\
L_+'&=\{j\mid\alpha_1+\alpha_2+\beta_j>1\}
&L_0'&=\{j\mid\alpha_1+\alpha_2+\beta_j=1\} \\
L_+^{(2)}&=L_+\backslash L_+^{(1)}
&L_0^{(2)}&=L_0\backslash L_0^{(1)},
\endalignedat \tag{\Lvarious} $$
we have
$$ \Afr_1=\smp{\Afr_{1,0}}{r_{1,0}}\oplus\bigoplus_{j\in L_+'}
\smdp\Cpx{\alpha_1+\alpha_2+\beta_j-1}{(p_1+p_2)\wedge q_j}\oplus
\bigoplus_{(i,j)\in L_+^{(1)}}\smdp\Cpx{\alpha_i+\beta_j-1}{p_i\wedge q_j},
$$
and there is a diffuse abelian subalgebra of the centralizer of
$\phi\restrict_{\Afr_{1,0}}$ containing $r_{1,0}(p_1+p_2)$.
By Proposition~\FreeProdDirectSum, $(p_1+p_2)\Afr(p_1+p_2)$ is freely
generated by $(p_1+p_2)\Afr_1(p_1+p_2)$ and $(\Cpx p_1+\Cpx p_2)$, so
$$ (p_1+p_2)\Afr(p_1+p_2)\cong
\left((p_1+p_2)\Afr_{1,0}(p_1+p_2)\oplus\bigoplus_{j\in L_+'}
\smdp\Cpx{\frac{\alpha_1+\alpha_2+\beta_j-1}{\alpha_1+\alpha_2}}
{(p_1+p_2)\wedge q_j}\right)
*\left(\smdp\Cpx{\frac{\alpha_1}{\alpha_1+\alpha_2}}{p_1}\oplus
\smdp\Cpx{\frac{\alpha_2}{\alpha_1+\alpha_2}}{p_2}\right). $$
Noting that $|L_+'|\le2$, we may use one of Lemma~\AzeroCCCC, Lemma~\AzeroCCC{}
or results from~\S\InPresOfDiff{} to show
$$ (p_1+p_2)\Afr(p_1+p_2)=
\smp{\Afr_{2,0}}{r_{2,0}}\oplus\bigoplus_{(i,j)\in L_+^{(2)}}
\smdp\Cpx{\alpha_i+\beta_j-1}{p_i\wedge q_j}. $$
Hence
$$ \Afr=\smp{\Afr_0}{r_0}\oplus\bigoplus_{(i,j)\in L_+}
\smdp\Cpx{\alpha_i+\beta_j-1}{p_i\wedge q_j}, \tag{\CdCCCrefpt} $$
where $r_0=r_{2,0}+r_{1,0}(1-p_1-p_2)$ and the linear span of the set
in~(\CaseIIAzeroCCCCdense) is dense in $\Afr_0$.

The inductive hypothesis (or Proposition~\TwoProj) yields
for every $(i,j)\in L_0^{(1)}$ a $*$--ho\-mo\-morph\-ism
$\pi_{(i,j)}^{(1)}:\Afr_{1,0}\to\Cpx$ such that
$\pi_{(i,j)}^{(1)}(r_{1,0}p_i)=1=\pi_{(i,j)}^{(1)}(r_{1,0}q_j)$, and for every
$(i,j)\in L_0^{(2)}$ a $*$--homomorphism
$\pi_{(i,j)}^{(2)}:\Afr_{2,0}\to\Cpx$ such that
$\pi_{(i,j)}^{(2)}(r_{2,0}p_i)=1=\pi_{(i,j)}^{(2)}(r_{2,0}q_j)$.
Looking at~(\CaseIIAzeroCCCCdense), one easily sees that each of these
$*$--homomorphisms can be uniquely extended to a $*$--homomorphism,
$\pi_{(i,j)}:\Afr_0\to\Cpx$, so that
$\pi_{(i,j)}(r_0p_i)=1=\pi_{(i,j)}(r_0q_j)$.

Let
$$ \align
\Afr_{1,00}&=\Afr_{1,0}\cap\bigcap_{(i,j)\in L_0^{(1)}}\ker\pi_{(i,j)}^{(1)},
 \\
\Afr_{2,00}&=\Afr_{2,0}\cap\bigcap_{(i,j)\in L_0^{(2)}}\ker\pi_{(i,j)}^{(2)}.
\endalign $$
Since $r_{1,0}(p_1+p_2)\in\Afr_{1,00}\cap\Afr_{2,00}$,
we see
from the denseness of the span of~(\CaseIIAzeroCCCCdense) in $\Afr_0$ that the
linear span of
$$ \aligned
\Afr_{2,0}&+\Afr_{2,00}\Afr_{1,00}(1-p_1-p_2)
 +(1-p_1-p_2)\Afr_{1,00}\Afr_{2,00}+\\
&+(1-p_1-p_2)\Afr_{1,00}\Afr_{2,00}\Afr_{1,00}(1-p_1-p_2)+
(1-p_1-p_2)\Afr_{1,0}(1-p_1-p_2)
\endaligned \tag{\CaseIICdCCCdense} $$
is dense in $\Afr_0$
and hence the linear span of
$$ \aligned
\Afr_{2,00}&+\Afr_{2,00}\Afr_{1,00}(1-p_1-p_2)
 +(1-p_1-p_2)\Afr_{1,00}\Afr_{2,00}+\\
&+(1-p_1-p_2)\Afr_{1,00}\Afr_{2,00}\Afr_{1,00}(1-p_1-p_2)+
(1-p_1-p_2)\Afr_{1,00}(1-p_1-p_2)
\endaligned $$
is dense in $\Afr_{00}$.
Note that $r_{1,0}(p_1+p_2)$ is full in $\Afr_{1,00}$.
Since $r_{1,0}(p_1+p_2)\in\Afr_{2,00}$, this implies that $\Afr_{2,00}$ is full
in $\Afr_{00}$.
Since $\Afr_{2,00}$ is simple, it follows from Proposition~\FullSimple{} that
$\Afr_{00}$ is simple.
(This also shows that $\Afr_0$ is simple when $L_0=\emptyset$.)

Let us now prove part~(ii).
If $i\in\{1,2\}$
then the linear span of
$$ \aligned
\bigg(\Afr_{2,0}\cap\bigcap
\Sb(i',j)\in L_0^{(2)}\\i'\ne i\endSb\ker_{(i',j)}^{(2)}\bigg)
&+\Afr_{2,00}\Afr_{1,00}(1-p_1-p_2)
 +(1-p_1-p_2)\Afr_{1,00}\Afr_{2,00}+\\ \vspace{2ex}
+(1-p_1-p_2)\Afr_{1,00}\Afr_{2,00}&\Afr_{1,00}(1-p_1-p_2)+
(1-p_1-p_2)\Afr_{1,00}(1-p_1-p_2)
\endaligned $$
is dense in~(\CdCCCdensei).
But $r_{2,0}p_i$ is full in
$$ \Afr_{2,0}\cap\bigcap\Sb(i',j)\in L_0^{(2)}\\i'\ne
i\endSb\ker_{(i',j)}^{(2)}. $$
Since this latter algebra contains $r_{1,0}(p_1+p_2)$, which is full in
$\Afr_{1,00}$ it then follows that $r_0p_i$ is full in the
algebra~(\CdCCCdensei).
Now take $i\in\{3,4,\ldots,n\}$.
For $j\in L_0'$ let $\pi_{(0,j)}^{(1)}:\Afr_{1,0}\to\Cpx$
be the $*$--homomorphism such that
$\pi_{(0,j)}^{(1)}(p_1+p_2)=1=\pi_{(0,j)}^{(1)}(q_j)$.
We have that $r_{1,0}p_i$ is full in
$$ \Afr_{1,0}\cap
\bigcap_{j\in L_0'}\ker\pi_{(0,j)}^{(1)}\cap
\bigcap\Sb(i',j)\in L_0^{(1)}\\i'\ne i\endSb
\ker\pi_{(i',j)}^{(1)}, $$
which in turn contains $(1-p_1-p_2)\Afr_{1,0}(p_1+p_2)$ and
$$ (1-p_1-p_2)\biggl(\Afr_{1,0}\cap
\bigcap\Sb(i',j)\in L_0^{(1)}\\i'\ne i\endSb
\ker\pi_{(i',j)}^{(1)}\biggr)(1-p_1-p_2). $$
But $(p_1+p_2)(\bigcap_{j\in L_0'}\ker\pi_{(0,j)}^{(1)})(p_1+p_2)$ meets
$\Afr_{2,00}$, which is simple.
Hence any ideal of the algebra~(\CdCCCdensej) which contains $r_0p_i$ must
also contain
$$ \aligned
\Afr_{2,00}
&+\Afr_{2,00}\Afr_{1,00}(1-p_1-p_2)
 +(1-p_1-p_2)\Afr_{1,00}\Afr_{2,00}+\\
&+(1-p_1-p_2)\Afr_{1,00}\Afr_{2,00}\Afr_{1,00}(1-p_1-p_2)+\\
&+(1-p_1-p_2)\bigg(\Afr_{1,0}\cap\bigcap\Sb(i',j)\in L_0^{(1)}\\i'\ne i\endSb
\ker\pi_{(i',j)}^{(1)}\bigg)(1-p_1-p_2),
\endaligned $$
which is dense in the algebra~(\CdCCCdensej).
This shows that $r_0p_i$ is full in~(\CdCCCdensej).
 
We now prove part~(iii).
We have that $r_{1,0}q_j$ is full in
$$ \Afr_{1,0}\cap\bigcap\Sb j'\in L_0'\\j'\ne j\endSb
\ker\pi_{(0,j')}^{(1)}\cap
\bigcap\Sb(i,j')\in L_0^{(1)}\\j'\ne j\endSb
\ker\pi_{(i,j')}^{(1)}, \tag{\Aonezerocapj} $$
which in turn contains $(1-p_1-p_2)\Afr_{1,0}(p_1+p_2)$ and
$$ (1-p_1-p_2)\biggl(\Afr_{1,0}\cap
\bigcap\Sb(i,j')\in L_0^{(1)}\\j'\ne j\endSb
\ker\pi_{(i,j')}^{(1)}\biggr)(1-p_1-p_2). $$
If $\exists i$ such that $(i,j)\in L_0^{(2)}$ then
$\alpha_1+\alpha_2+\beta_j>1$ so $(p_1+p_2)\wedge q_j\ne0$ and
$r_{2,0}((p_1+p_2)\wedge q_j)$ is full in
$$ \Afr_{2,0}\cap\bigcap\Sb(i,j')\in L_0^{(2)}\\j'\ne j\endSb
\ker\pi_{(i,j')}^{(2)}. \tag{\Atwozerocapj} $$
Hence any ideal of the algebra~(\CdCCCdensej) that contains $r_0q_j$ must
contain
$$ \aligned
\bigg(\Afr_{2,0}\cap\bigcap
\Sb(i,j')\in L_0^{(2)}\\j'\ne j\endSb&\ker_{(i,j')}^{(2)}\bigg)
+\Afr_{2,00}\Afr_{1,00}(1-p_1-p_2)
+(1-p_1-p_2)\Afr_{1,00}\Afr_{2,00}+\\ \vspace{1.5ex}
&+(1-p_1-p_2)\Afr_{1,00}\Afr_{2,00}\Afr_{1,00}(1-p_1-p_2)+\\ \vspace{1ex}
&+(1-p_1-p_2)\bigg(\Afr_{1,0}\cap\bigcap\Sb(i,j')\in L_0^{(1)}\\j'\ne j\endSb
\ker_{(i,j')}^{(1)}\bigg)(1-p_1-p_2),
\endaligned $$
whose span is dense in~(\CdCCCdensej).
On the other hand, if there is no $i$ such that $(i,j)\in L_0^{(2)}$ then the
algebra~(\Atwozerocapj) is $\Afr_{2,00}$, which is simple.
By a dimension argument, the algebra~(\Aonezerocapj) meets $\Afr_{2,00}$.
Hence any ideal of the algebra~(\CdCCCdensej) that contains $r_0q_j$ must
contain
$$ \aligned
\Afr_{2,00}&+\Afr_{2,00}\Afr_{1,00}(1-p_1-p_2)
+(1-p_1-p_2)\Afr_{1,00}\Afr_{2,00}+\\
&+(1-p_1-p_2)\Afr_{1,00}\Afr_{2,00}\Afr_{1,00}(1-p_1-p_2)+\\
&+(1-p_1-p_2)\bigg(\Afr_{1,0}\cap\bigcap\Sb(i,j')\in L_0^{(1)}\\j'\ne j\endSb
\ker_{(i,j')}^{(1)}\bigg)(1-p_1-p_2),
\endaligned $$
whose span is dense in~(\CdCCCdensej).
Thus $r_0q_j$ is full in~(\CdCCCdensej).

The required results about stable rank and uniqueness of the trace
follow as in previous lemmas from the fact that $\Afr_{2,00}$ is full in
$\Afr_{00}$.
\QED

\proclaim{Lemma \AzeroCdCCC}
Let $n\in\Naturals\cup\{0\}$ and let
$$ (\Afr,\phi)=(A_0\oplus\smdp\Cpx{\alpha_1}{p_1}\oplus\cdots\oplus
\smdp\Cpx{\alpha_n}{p_n})*
(\smdp\Cpx{\beta_1}{q_1}\oplus\smdp\Cpx{\beta_2}{q_2}), $$
where the centralizer of $\phi\restrict_{A_0}$ has a diffuse abelian
subalgebra.
Let $L_+$ and $L_0$ be as in~(\LplusLzeron).
Then
$$ \Afr=\smp{\Afr_0}{r_0}\oplus
\bigoplus_{(i,j)\in L_+}\smdp\Cpx{\alpha_i+\beta_j-1}{p_i\wedge q_j}, 
\tag{\AzeroCdCCCAfr} $$
where the centralizer of $\phi\restrict_{\Afr_0}$ has a diffuse abelian
subalgebra which contains $r_0p_1$ and a diffuse abelian subalgebra
which contains $r_0q_1$.

If $\phi\restrict_{A_0}$ is a trace then the stable rank of $\Afr$ is~1.

If $L_0$ is empty then $\Afr_0$ is simple.
If, in addition, $\phi\restrict_{A_0}$ is a trace then
$\phi(r_0)^{-1}\phi\restrict_{\Afr_0}$
is the unique tracial state
on $\Afr_0$ and if $\phi\restrict_{A_0}$ is not a trace then $\Afr_0$ has no
tracial states.

If $L_0$ is not empty, then for every $(i,j)\in L_0$ there is a
$*$--homomorphism $\pi_{(i,j)}:\Afr_0\to\Cpx$ such that
$\pi_{(i,j)}(r_0p_i)=1=\pi_{(i,j)}(r_0q_j)$.
Then
\roster
\item"(i)"
$$ \Afr_{00}\eqdef\dsize\bigcap_{(i,j)\in L_0}\ker\pi_{(i,j)} $$
is simple.
If $\phi\restrict_{A_0}$ is a trace then
$\phi(r_0)^{-1}\phi\restrict_{\Afr_{00}}$
is the unique tracial state
on $\Afr_0$ and if $\phi\restrict_{A_0}$ is not a trace then $\Afr_{00}$ has no
tracial states.
\item"(ii)" For each $i\in\{1,2,\ldots,n\}$, $r_0p_i$ is full in
$$ \Afr_0\cap\bigcap\Sb(i',j)\in L_0\\i'\ne i\endSb \ker\pi_{(i',j)}.
\tag{\AzeroCdCCCdensei} $$
\vskip2ex
\item"(iii)" For each $j\in\{1,2\}$, $r_0q_j$ is full in
$$ \Afr_0\cap\bigcap\Sb(i,j')\in L_0\\j'\ne j\endSb\ker\pi_{(i,j')}.
\tag{\AzeroCdCCCdensej} $$
\endroster
\endproclaim
\demo{Proof}
The cases $n=0,1,2$ have been already proved.
Let $p_0=1-\sum_1^np_j$ and let $\Afr_1$ be the C$^*$--subalgebra of $\Afr$
generated by $(\Cpx p_0+\Cpx p_1+\cdots+\Cpx p_n)\cup(\Cpx q_1+\Cpx q_2)$, so
that
$$ (\Afr_1,\phi\restrict_{\Afr_1})=(\smdp\Cpx{\alpha_0}{p_0}
\oplus\smdp\Cpx{\alpha_1}{p_1}\oplus\cdots\oplus\smdp\Cpx{\alpha_n}{p_n})*
(\smdp\Cpx{\beta_1}{q_1}\oplus\smdp\Cpx{\beta_2}{q_2}). $$
Let
$$ \align
L_+^{(0)}&=\{(0,j)\mid\alpha_0+\beta_j>1\}\\
L_0^{(0)}&=\{(0,j)\mid\alpha_0+\beta_j=1\}
\endalign $$
We use Lemma~\CdCCC{} to find
$$ \Afr_1=\smp{\Afr_{1,0}}{r_{1,0}}\oplus\bigoplus_{(i,j)\in L_+\cup L_+^{(0)}}
\smdp\Cpx{\alpha_i+\beta_j-1}{p_i\wedge q_j}. $$
Then by Proposition~\FreeProdDirectSum, $p_0\Afr_1p_0$ and $A_0$ freely
generate $p_0\Afr p_0$.
Hence by Proposition~\FPHaarUnitaryDixmier, $p_0\Afr p_0$ is simple.
So~(\AzeroCdCCCAfr) holds where $r_0=p_0+(1-p_0)r_{1,0}$ and where the linear
span of
$$ p_0\Afr p_0+p_0\Afr p_0\Afr_{1,0}(1-p_0)+(1-p_0)\Afr_{1,0}p_0\Afr p_0+
(1-p_0)\Afr_{1,0}p_0\Afr p_0\Afr_{1,0}(1-p_0)+(1-p_0)\Afr_{1,0}(1-p_0)
\tag{\AzeroCdCCCdense} $$
is dense in $\Afr_0$.
The application of Lemma~\CdCCC{} gives for each $(i,j)\in L_0\cup L_0^{(0)}$ a
$*$--ho\-mo\-morph\-ism $\pi_{(i,j)}^{(1)}:\Afr_{1,0}\to\Cpx$ such that
$\pi_{(i,j)}^{(1)}(r_{1,0}p_i)=1=\pi_{(i,j)}^{(1)}(r_{1,0}q_j)$
and then $r_{1,0}p_0$ is full in
$$ \Afr_{1,0}\cap\bigcap_{(i,j)\in L_0}\ker\pi_{(i,j)}^{(1)}.
\tag{\AzeroCdCCCAfronefull} $$
For each $(i,j)\in L_0$, $\pi_{(i,j)}^{(1)}$ extends to a $*$--homomorphism
$\pi_{(i,j)}:\Afr_0\to\Cpx$ whose kernel is densely spanned by
$$ \align
p_0\Afr p_0&+p_0\Afr p_0\Afr_{1,0}(1-p_0)+(1-p_0)\Afr_{1,0}p_0\Afr p_0+ \\
&+(1-p_0)\Afr_{1,0}p_0\Afr p_0\Afr_{1,0}(1-p_0)
+(1-p_0)(\ker\pi_{(i,j)}^{(1)})(1-p_0).
\endalign $$
Since the algebra~(\AzeroCdCCCAfronefull) contains both
$(1-p_0)\big(\Afr_{1,0}\cap\bigcap_{(i,j)\in L_0}
\ker\pi_{(i,j)}^{(1)}(1-p_0)\big)$
and
$(1-p_0)\Afr_{1,0}p_0$, from the denseness of~(\AzeroCdCCCdense) in
$\Afr_0$ we see that $p_0$ is full in $\Afr_{00}$.
Since $p_0\Afr p_0$ is simple, by Proposition~\FullSimple{} this implies that
$\Afr_{00}$ is simple, (hence when $L_0$ is empty, $\Afr_0$ is simple).

From the facts, for $(i,j)\in L_0$, that $r_{1,0}p_i$ is full in
$$ \Afr_{1,0}\cap\bigcap_{(0,j)\in L_0^{(0)}}\ker\pi_{(0,j)}^{(1)}\cap
\bigcap\Sb(i'.j)\in L_0\\i'\ne i\endSb\ker\pi_{(i',j)}^{(1)}, $$
which by a dimension argument contains a nonzero element of $p_0\Afr_{1,0}p_0$,
and that every positive element of $p_0\Afr p_0$ is full in $\Afr_{00}$,
we obtain that $r_0p_i$ is full in~(\AzeroCdCCCdensei), proving~(ii).

Because $r_{1,0}q_j$ is full in
$$ \Afr_{1,0}\cap\bigcap\Sb(0,j')\in L_0^{(0)}\\j'\ne j\endSb
\ker\pi_{(0,j)}^{(1)}\cap
\bigcap\Sb(i.j')\in L_0\\j'\ne j\endSb\ker\pi_{(i,j')}^{(1)}, $$
which meets $p_0\Afr_{1,0}p_0$, we similarly obtain that $r_0q_j$ is full
in~(\AzeroCdCCCdensej), proving~(iii).

If $\phi\restrict_{A_0}$ is a trace
then, by Propositions~\FPHaarUnitaryDixmier{} and~\FPHaarUnitarySRone,
$p_0\Afr p_0$ has unique tracial state and stable rank~1.
Since $p_0\Afr p_0$ is full in $\Afr_{00}$, by Proposition~\AheredBsrone{} the
same hold for $\Afr_{00}$ (which is just $\Afr_0$ when $L_0$ is empty)
and $\phi(r_0)^{-1}\phi\restrict_{\Afr_{00}}$ is then seen to be the unique
tracial state.

If $\phi\restrict_{A_0}$ is not a trace then by
Proposition~\FPHaarUnitaryDixmier, $p_0\Afr p_0$ has no tracial state, so
neither does $\Afr_{00}$ have a tracial state.
\QED

The following proposition proves Theorem~\FPfdAbelian.

\proclaim{Proposition \CdCCdC}
Let
$$ (\Afr,\phi)=(\smdp\Cpx{\alpha_1}{p_1}\oplus\cdots\oplus
\smdp\Cpx{\alpha_n}{p_n})*
(\smdp\Cpx{\beta_1}{q_1}\oplus\cdots\oplus
\smdp\Cpx{\beta_m}{q_m}), $$
where $n\ge2$ and $m\ge3$.
Let
$$ \aligned
L_+&=\{(i,j)\mid\alpha_i+\beta_j>1\}\\
L_0&=\{(i,j)\mid\alpha_i+\beta_j=1\}.
\endaligned $$
Then
$$ \Afr=\smp{\Afr_0}{r_0}\oplus
\bigoplus_{(i,j)\in L_+}\smdp\Cpx{\alpha_i+\beta_j-1}{p_i\wedge q_j} $$
where the centralizer of $\phi\restrict_{\Afr_0}$ has a diffuse abelian
subalgebra containing $r_0p_1$.

If $L_0$ is empty then $\Afr_0$ is simple, nonunital and
$\phi(r_0)^{-1}\phi\restrict_{\Afr_0}$ is the unique tracial state on
$\Afr_0$.

If $L_0$ is not empty, then for every $(i,j)\in L_0$ there is a
$*$--homomorphism $\pi_{(i,j)}:\Afr_0\to\Cpx$ such that
$\pi_{(i,j)}(r_0p_i)=1=\pi_{(i,j)}(r_0q_j)$.
Then
\roster
\item"(i)"
$$ \Afr_{00}\eqdef\dsize\bigcap_{(i,j)\in L_0}\ker\pi_{(i,j)} $$
is
simple and nonunital and $\phi(r_0)^{-1}\phi\restrict_{\Afr_{00}}$ is the unique
tracial state on $\Afr_{00}$.
\vskip2ex
\item"(ii)" For each $i\in\{1,2,\ldots,n\}$, $r_0p_i$ is full in
$\dsize\Afr_0\cap\bigcap\Sb(i',j)\in L_0\\i'\ne i\endSb \ker\pi_{(i',j)}$.
\endroster
\endproclaim
\demo{Proof}
We proceed by induction on $\min(n,m)$.
If $\min(n,m)=2$ then Lemma~\CdCCC{} applies to prove the desired results.
If $\min(n,m)\ge3$ then take $n\le m$ and let
$\Afr_1$ be the C$^*$--subalgebra of $\Afr$ generated by
$(\Cpx(p_1+p_2)+\Cpx p_3+\cdots\Cpx p_n)\cup(\Cpx q_1+\cdots+\Cpx q_m)$.
Thus
$$ (\Afr_1,\phi\restrict_{\Afr_1})\cong
(\smdp\Cpx{\alpha_1+\alpha_2}{p_1+p_2}\oplus\smdp\Cpx{\alpha_3}{p_3}
\oplus\cdots\oplus\smdp\Cpx{\alpha_n}{p_n})
*(\smdp\Cpx{\beta_1}{q_1}\oplus\cdots\oplus\smdp\Cpx{\beta_m}{q_m}). $$
By inductive hypothesis, letting $L_+$, $L_0$, $L_+^{(1)}$, $L_0^{(1)}$,
$L_+'$, $L_0'$, $L_+^{(2)}$ and $L_0^{(2)}$ be as in~(\Lvarious),
we have
$$ \Afr_1=\smp{\Afr_{1,0}}{r_{1,0}}\oplus\bigoplus_{j\in L_+'}
\smdp\Cpx{\alpha_1+\alpha_2+\beta_j-1}{(p_1+p_2)\wedge q_j}\oplus
\bigoplus_{(i,j)\in L_+^{(1)}}\smdp\Cpx{\alpha_i+\beta_j-1}{p_i\wedge q_j},
$$
and there is a diffuse abelian subalgebra of the centralizer of
$\phi\restrict_{\Afr_{1,0}}$ containing $r_{1,0}(p_1+p_2)$.
By Proposition~\FreeProdDirectSum, $(p_1+p_2)\Afr(p_1+p_2)$ is freely
generated by $(p_1+p_2)\Afr_1(p_1+p_2)$ and $(\Cpx p_1+\Cpx p_2)$, so
$$ (p_1+p_2)\Afr(p_1+p_2)\cong
\left((p_1+p_2)\Afr_{1,0}(p_1+p_2)\oplus\bigoplus_{j\in L_+'}
\smdp\Cpx{\frac{\alpha_1+\alpha_2+\beta_j-1}{\alpha_1+\alpha_2}}
{(p_1+p_2)\wedge q_j}\right)
*\left(\smdp\Cpx{\frac{\alpha_1}{\alpha_1+\alpha_2}}{p_1}\oplus
\smdp\Cpx{\frac{\alpha_2}{\alpha_1+\alpha_2}}{p_2}\right). $$
Applying Lemma~\AzeroCdCCC{} yields
$$ (p_1+p_2)\Afr(p_1+p_2)=
\smp{\Afr_{2,0}}{r_{2,0}}\oplus\bigoplus_{(i,j)\in L_+^{(2)}}
\smdp\Cpx{\alpha_i+\beta_j-1}{p_i\wedge q_j}. $$
Hence
$$ \Afr=\smp{\Afr_0}{r_0}\oplus\bigoplus_{(i,j)\in L_+}
\smdp\Cpx{\alpha_i+\beta_j-1}{p_i\wedge q_j}, $$
where $r_0=r_{2,0}+r_{1,0}(1-p_1-p_2)$ and the linear span of the set
in~(\CaseIIAzeroCCCCdense) is dense in $\Afr_0$.

Now the proof of this proposition follows verbatim the proof of Lemma~\CdCCC{}
after equation~(\CdCCCrefpt), except we must also show that when $L_0$ is
nonempty
then $\Afr_{00}$ is nonunital.
Suppose for contradiction that $\Afr_{00}$ is unital with identity $e$.
Then $e\ne r_0$ but $ex=exe=xe$ for every $x\in\Afr_0$.
Thus $e$ is in the center of $\Afr_0$.
But by the results of~\cite{\DykemaZZFreeDim}, the strong--operator closure of
$\Afr_0$ (in the GNS representation associated to $\phi\restrict_{\Afr_0}$) is
a II$_1$--factor.
This gives a contradiction, because $e$ must be in the center of this von
Neumann algebra.
\QED

\proclaim{Remark \NonUnitalToo}\rm
By the same proof, one shows that for every nonempty subset, $F$, of $L_0$, the
ideal $\bigcap_{(i,j)\in F}\ker\pi_{(i,j)}$ of $\Afr_0$ is nonunital.
\endproclaim

The following lemma can be proved from Proposition~\CdCCdC{} using
Proposition~\FreeProdDirectSum{} in the same way
that Lemma~\AzeroCdCCC{} was proved from Lemma~\CdCCC.
\proclaim{Lemma \AzeroCdCCdC}
Let
$$ (\Afr,\phi)=(A_0\oplus\smdp\Cpx{\alpha_1}{p_1}\oplus\cdots\oplus
\smdp\Cpx{\alpha_n}{p_n})*
(\smdp\Cpx{\beta_1}{q_1}\oplus\cdots\oplus
\smdp\Cpx{\beta_m}{q_m}), \tag{\eqAzeroCdCCdC} $$
where the centralizer $\phi\restrict_{A_0}$ has a diffuse abelian subalgebra
and where $n\ge1$, $m\ge2$.
Or, let
$$ (\Afr,\phi)=(A_0\oplus\smdp\Cpx{\alpha_1}{p_1}\oplus\cdots\oplus
\smdp\Cpx{\alpha_n}{p_n})*
(B_0\oplus\smdp\Cpx{\beta_1}{q_1}\oplus\cdots\oplus
\smdp\Cpx{\beta_m}{q_m}), \tag{\eqAzeroCdCBzeroCdC} $$
where the centralizer of each of $\phi\restrict_{A_0}$ and
$\phi\restrict_{B_0}$ has a diffuse abelian subalgebra and where $n\ge1$,
$m\ge1$.

Then
$$ \Afr=\smp{\Afr_0}{r_0}\oplus
\bigoplus_{(i,j)\in L_+}\smdp\Cpx{\alpha_i+\beta_j-1}{p_i\wedge q_j} $$
where the centralizer of $\phi\restrict_{\Afr_0}$ has a diffuse abelian
subalgebra containing $r_0p_1$ and another containing $r_0q_1$.

If $L_0$ is empty then $\Afr_0$ is simple.

If $L_0$ is not empty, then for every $(i,j)\in L_0$ there is a
$*$--homomorphism $\pi_{(i,j)}:\Afr_0\to\Cpx$ such that
$\pi_{(i,j)}(r_0p_i)=1=\pi_{(i,j)}(r_0q_j)$.
Then
\roster
\item"(i)"
$$\Afr_{00}\eqdef\dsize\bigcap_{(i,j)\in L_0}\ker\pi_{(i,j)} $$
is
simple and nonunital and $\phi(r_0)^{-1}\phi\restrict_{\Afr_{00}}$ is the unique
tracial state on $\Afr_{00}$.
\vskip2ex
\item"(ii)"  For each $i\in\{1,2,\ldots,n\}$, $r_0p_i$ is full in
$\dsize\Afr_0\cap\bigcap\Sb(i',j)\in L_0\\i'\ne i\endSb \ker\pi_{(i',j)}$.
\item"(iii)" For each $j\in\{1,2,\ldots,m\}$, $r_0p_j$ is full in
$\dsize\Afr_0\cap\bigcap\Sb(i,j')\in L_0\\j'\ne j\endSb \ker\pi_{(i,j')}$.
\endroster
Moreover, if $\phi_{A_0}$ is a trace (and if, in the case
of~(\eqAzeroCdCBzeroCdC),
$\phi\restrict_{B_0}$ is a trace) then $\Afr$ has stable rank~1 and
$\phi(r_0)^{-1}\phi$ is the unique
tracial state on $\Afr_0$ when $L_0$ is empty or $\Afr_{00}$ when $L_0$ is
nonempty.
Otherwise, $\Afr_0$, respectively $\Afr_{00}$, has no tracial state.
\endproclaim

As promised in~\S\StmntMainResults, a result similar to~\CdCCdC{} holds for
free products of more than two finite dimensional abelian C$^*$--algebras.
\proclaim{Proposition \FinManyFDAbel}
Let $N\in\Naturals$, $N\ge3$ and for
each $\iota\in\{1,\ldots,N\}$ take a finite dimensional
abelian C$^*$--algebra and faithful tracial state,
$$ (A_\iota,\tau_\iota)=
\smdp\Cpx{\alpha_{\iota,1}}{p_{\iota,1}}\oplus
\smdp\Cpx{\alpha_{\iota,2}}{p_{\iota,2}}\oplus\cdots\oplus
\smdp\Cpx{\alpha_{\iota,n(\iota)}}{p_{\iota,n(\iota)}}, \tag{\Aiota} $$
where $n(\iota)\in\Naturals$, $n(\iota)\ge2$.
Let
$$ (\Afr,\tau)=\operatornamewithlimits*_{\iota=1}^N(A_\iota,\tau_\iota) $$
and let
$$ \aligned
L_+&=\left\{\bigl(j(\iota)\bigr)_{\iota=1}^N\;\bigg|\;
\sum_{\iota=1}^N(1-\alpha_{\iota,j(\iota)})<1\right\} \\ \vspace{1ex}
L_0&=\left\{\bigl(j(\iota)\bigr)_{\iota=1}^N\;\bigg|\;
\sum_{\iota=1}^N(1-\alpha_{\iota,j(\iota)})=1\right\}.
\endaligned \tag{\LpluszeroMTTFDA} $$
Then
$$ \Afr=\smp{\Afr_0}{r_0}\oplus
\bigoplus_{j\in L_+}\smdp\Cpx{\gamma_j}{r_j} $$
where for $j=(j(\iota))_{\iota=1}^N\in L_+$,
$\gamma_j=1-\sum_{\iota=1}^N(1-\alpha_{\iota,j(\iota)})$ and
$r_j=\bigwedge_{\iota=1}^Np_{\iota,j(\iota)}$.
For each $1\le\iota\le N$ and each $1\le k\le n(\iota)$, $\Afr_0$ has a diffuse
abelian subalgebra which contains $r_0p_{\iota,k}$.
Moreover, $\Afr$ has stable rank~1.

If $L_0$ is empty then $\Afr_0$ is simple and
$\phi(r_0)^{-1}\phi\restrict_{\Afr_0}$ is the unique tracial state on $\Afr_0$.

If $L_0$ is nonempty then for every $j=(j(\iota))_{\iota=1}^N\in L_0$ there is
a $*$--homomorphism $\pi_j:\Afr_0\to\Cpx$ such that $\forall1\le\iota\le N$
$\pi_j(p_{\iota,j(\iota)})=1$.
Then
\roster
\item"(i)"
$$ \Afr_{00}\eqdef\bigcap_{j\in L_0}\ker\pi_j $$
is simple and nonunital and
$\phi(r_0)^{-1}\phi\restrict_{\Afr_{00}}$ is the unique tracial state on
$\Afr_{00}$.
\item "(ii)" For each $1\le\iota\le N$ and each $1\le k\le n(\iota)$,
$r_0p_{\iota,k}$ is full in
$$ \Afr_0\cap\bigcap\Sb j\in L_0\\j(\iota)\ne k\endSb\ker\pi_j. $$
\endroster
\endproclaim
\demo{Proof}
The proof is by induction on $N$, and the case $N=3$ and inductive step are
proved simultaneously.
Let $\Afr_{N-1}$ be the C$^*$--subalgebra of $\Afr$ generated by
$\bigcup_{\iota=1}^{N-1}A_\iota$.
Use the inductive hypothesis (or, when $N=3$, Proposition~\CdCCdC{} or
Proposition~\TwoProj) to find
$$ (\Afr_{N-1},\tau\restrict_{\Afr_{N-1}})\cong
{\operatornamewithlimits{\ast}_{\iota=1}^{N-1}}(A_\iota,\tau_\iota). $$
Then
$$ (\Afr,\tau)\cong(\Afr_{N-1},\tau\restrict_{\Afr_{N-1}})*(A_N,\tau_N), $$
and one uses Lemma~\AzeroCdCCdC{} to find $\Afr$.
\QED

One can generalize this to the free product of finitely many (say $N$) algebras
of the
form $A_0\oplus\Cpx\oplus\cdots\oplus\Cpx$, either, as
Proposition~\FinManyFDAbel{} was proved, by using
Lemma~\AzeroCdCCdC{} and induction on $N$, or, as Lemma~\AzeroCdCCC{}
was proved, by applying Proposition~\FinManyFDAbel{} and
Lemma~\FPHaarUnitaryDixmier.
One obtains the following result.
\proclaim{Proposition \FinManyAzeroCdC}
Let $N\in\Naturals$, $N\ge3$ and for
each $\iota\in\{1,\ldots,N\}$ let $(A_\iota,\phi_\iota)$ be either a finite
dimensional abelian algebra with a faithful state as in~(\Aiota) or
$$ (A_\iota,\phi_\iota)=A_{\iota,0}\oplus
\smdp\Cpx{\alpha_{\iota,1}}{p_{\iota,1}}\oplus
\smdp\Cpx{\alpha_{\iota,2}}{p_{\iota,2}}\oplus\cdots\oplus
\smdp\Cpx{\alpha_{\iota,n(\iota)}}{p_{\iota,n(\iota)}}, \tag{\AzeroCdC} $$
where $n(\iota)\in\Naturals$, $n(\iota)\ge1$ and where
the centralizer of $\phi_\iota\restrict_{A_{\iota,0}}$ has a diffuse abelian
subalgebra.
Let
$$ (\Afr,\tau)=\operatornamewithlimits*_{\iota=1}^N(A_\iota,\tau_\iota) $$
and let $L_+$ and $L_0$ be as in~(\LpluszeroMTTFDA).
Then
$$ \Afr=\smp{\Afr_0}{r_0}\oplus
\bigoplus_{j\in L_+}\smdp\Cpx{\gamma_j}{r_j} $$
where for $j=(j(\iota))_{\iota=1}^N\in L_+$,
$\gamma_j=1-\sum_{\iota=1}^N(1-\alpha_{\iota,j(\iota)})$ and
$r_j=\bigwedge_{\iota=1}^Np_{\iota,j(\iota)}$.
If each $\phi_\iota\restrict_{A_{\iota,0}}$ is a trace then $\Afr$ has stable
rank~1.

If $L_0$ is empty then $\Afr_0$ is simple and if each
$\phi_\iota\restrict_{A_{\iota,0}}$ is a trace then
$\phi(r_0)^{-1}\phi\restrict_{\Afr_0}$ is the unique tracial state on $\Afr_0$.
If some $\phi_\iota\restrict_{A_{\iota,0}}$ is not a trace then $\Afr_0$ has no
tracial states.

If $L_0$ is nonempty then for every $j=(j(\iota))_{\iota=1}^N\in L_0$ there is
a $*$--homomorphism $\pi_j:\Afr_0\to\Cpx$ such that $\forall1\le\iota\le N$
$\pi_j(p_{\iota,j(\iota)})=1$.
Then
$$ \Afr_{00}\eqdef\bigcap_{j\in L_0}\ker\pi_j $$
is simple and nonunital and if each
$\phi_\iota\restrict_{A_{\iota,0}}$ is a trace then
$\phi(r_0)^{-1}\phi\restrict_{\Afr_{00}}$ is the unique tracial state on
$\Afr_{00}$.
If some $\phi_\iota\restrict_{A_{\iota,0}}$ is not a trace then $\Afr_{00}$ has
no tracial states.
\endproclaim

Restricting the above proposition to the case when each $A_{\iota,0}$ is
abelian, we obtain the
following result about the free product of finitely many abelian
C$^*$--algebras.
\proclaim{Corollary \FinManyCofXAtoms}
Let $N\in\Naturals$, $N\ge2$ and for each
$1\le\iota\le N$ consider the abelian C$^*$--algebra and state
$(C(X_\iota),\int\cdot\dif\mu_\iota)$, where $\mu_\iota$ is a regular, Borel
probability measure on $X_\iota$ whose support is all of $X_\iota$ and having
at most finitely many atoms, each of which is an isolated point of $X_\iota$.
Let
$$ (\Afr,\tau)=\operatornamewithlimits*_{\iota=1}^N
(C(X_\iota),\tint\cdot\dif\mu_\iota). $$
Let
$$ \align
L_+&=\left\{x=(x_\iota)_{\iota=1}^N\bigg|x_\iota\in X_\iota,\,
\sum_{\iota=1}^N(1-\mu_\iota(\{x_\iota\}))<1\right\}\\ \vspace{2ex}
L_0&=\left\{x=(x_\iota)_{\iota=1}^N\bigg|x_\iota\in X_\iota,\,
\sum_{\iota=1}^N(1-\mu_\iota(\{x_\iota\}))=1\right\}.
\endalign $$
Assume that no $X_\iota$ is a one--point space and also
exclude the case when $N=2$ and $X_1$ and $X_2$ are both two--point spaces.
Then $\Afr$ has stable rank~1 and
$$ \Afr=\smp{\Afr_0}{r_0}\oplus\bigoplus_{x\in L_+}\smdp\Cpx{\alpha_x}{r_x}
\tag{\AfrCofXIsol} $$
where $\alpha_x=1-\sum_{\iota=1}^N(1-\mu_\iota(\{x_\iota\}))$.
If $L_0$ is empty the $\Afr_0$ is simple and has unique tracial state
$\tau(r_0)^{-1}\tau\restrict_{\Afr_0}$.

If $L_0$ is nonempty, then there are distinct, surjective  $*$--homomorphisms,
$\pi_x:\Afr_0\to\Cpx$, $(x\in L_0)$, such that
$$ \Afr_{00}\eqdef\bigcap_{x\in L_0}\ker\pi_x $$
is simple and nonunital and has unique tracial state
$\tau(r_0)^{-1}\tau\restrict_{\Afr_{00}}$.
\endproclaim

\vskip3ex
\noindent{\bf\S\CofXS.  More general abelian C$^*$--algebras.}
\vskip3ex

In this section we will investigate free products,
$$ (\Afr,\tau)=\freeprod{\iota=1}N(C(X_\iota),\tint\cdot\dif\mu_\iota), $$
of abelian
C$^*$--algebras and states, $(C(X_\iota),\tint\cdot\dif\mu_\iota)$, each of
which can be written as an inductive
limit of the abelian algebras and states considered in
Corollary~\FinManyCofXAtoms.
The criterion for simplicity of the free product of such abelian algebras is
the same as for finite dimensional abelian algebras, namely, $\Afr$ is simple
if and only if there are not atoms $x_\iota\in X_\iota$ of $\mu_\iota$
$(1\le\iota\le N)$ such that $\sum_1^N(1-\mu(\{x_\iota\}))\le1$.
However, in the case when there are atoms $x_\iota$ satisfying
$\sum_1^N(1-\mu(\{x_\iota\}))<1$, we don't always get a corresponding copy of
$\Cpx$ as a direct summand of $\Afr$.
In fact
we get such a direct summand, i.e\. a minimal and central projection 
in $\Afr$, like
$r_x$ in~(\AfrCofXIsol), corresponding to these atoms,
if and only if each $x_\iota$ is an
isolated point of $X_\iota$.

\proclaim{Definition \InvLimIsolAtomsDef}\rm
Let $X$ be a compact, Hausdorff, topological space and let $\mu$ be a regular,
Borel, probability measure on $X$.
We say $(X,\mu)$ is {\it an inverse limit of spaces and measures with isolated
atoms} if $X$ is an inverse limit, $X=\varprojlim(X_n,\kappa_n)$,
of compact Hausdorff spaces, $X_n$, and surjective, continuous maps,
$\kappa_n:X_{n+1}\to X_n$, $(n\in\Naturals)$, such that letting $\lambda_n:X\to
X_n$ be the
resulting canonical, surjective maps and letting $\mu_n=(\lambda_n)_*(\mu)$ be
the push--forward measures, each $\mu_n$ has at most finitely many atoms and
each atom of $\mu_n$ is an isolated point of $X_n$ and, moreover, if $x\in X_n$
and if $\kappa_n^{-1}(\{x\})$ has more than one point, then $x$ is an atom of
$\mu_n$.
\endproclaim

\proclaim{Examples \InvLimIsolAtomsEx}\rm
In each of the following cases, $(X,\mu)$ is an inverse
limit of spaces and measures with isolated atoms.
(One can easily cook up more intricate examples as well.)
\roster
\item"(i)" $\mu$ has no atoms
\item"(ii)" all the atoms of $\mu$ are isolated points of $X$
\item"(iii)" $X$ is separable and totally disconnected
\item"(iv)" $X=\{0\}\cup\bigcup_{n=1}^\infty[\frac1{2n+1},\frac1{2n}]$ with the
relative topology from $\Reals$, and $\mu$ has a unique atom at $0$.
\endroster
\endproclaim

\proclaim{Proposition \InvLimIsolAtomsFP}
Let $N\in\Naturals$, $N\ge2$, and for each $1\le\iota\le N$ let
$(X_\iota,\mu_\iota)$ be a compact Hausdorff space and a regular, Borel,
probability measures, each of which is an inverse limit of spaces and measures
with isolated atoms.
Assume each $X_\iota$ has more than one point and each $\mu_\iota$ has support
equal to all of $X_\iota$.
Exclude the case when $N=2$ and $X_1$ and $X_2$ are both two--point spaces.
Let
$$ (\Afr,\tau)=\freeprod{\iota=1}N(C(X_\iota),\tint\cdot\dif\mu_\iota). $$
Let
$$ \align
I_+&=\left\{x=(x_\iota)_{\iota=1}^N\in\prod_{\iota=1}^NX_\iota
\bigg|\sum_{\iota=1}^N(1-\mu_\iota(\{x_\iota\}))<1,\;
\text{ each $x_\iota$ is isolated in $X_\iota$}\right\} \\
L&=\left\{x=(x_\iota)_{\iota=1}^N\in\prod_{\iota=1}^NX_\iota
\bigg|\sum_{\iota=1}^N(1-\mu_\iota(\{x_\iota\}))\le1\right\}
\bigg\backslash I_+.
\endalign $$
Then $\Afr$ has stable rank~1 and
$$ \Afr=\smp{\Afr_0}{r_0}\oplus\bigoplus_{x\in I_+}\smdp\Cpx{\alpha_x}{r_x},
\tag{\AfrIndLim} $$
where $\alpha_x=1-\sum_{\iota=1}^N(1-\mu_\iota(\{x_\iota\}))$ and where
$r_x=\bigwedge_{\iota=1}^Np_{x_\iota}$, with $p_{x_\iota}\in C(X_\iota)$ the
characteristic function of $\{x_\iota\}$.
If $L$ is empty then $\Afr_0$ is simple and has unique tracial state
$\tau(r_0)^{-1}\tau\restrict_{\Afr_0}$.

If $L$ is nonempty then for every $x\in L$ there is a $*$--homomorphism
$\pi_x:\Afr_0\to\Cpx$ such that whenever $f\in C(X_\iota)$ we
have
$\pi_x(f)=f(x_\iota)$.
Moreover,
$$ \Afr_{00}\eqdef\bigcap_{x\in L}\ker\pi_x $$
is simple and nonunital, and has unique tracial state.
Finally, for each nonempty subset, $F\subseteq L$, the ideal
$\bigcap_{x\in F}\ker\pi_x$ of $\Afr_0$ is nonunital.
\endproclaim
\demo{Proof}
Let $X_\iota=\varprojlim(X_{\iota,n},\kappa_{\iota,n})$ be an inverse limit
with properties as in Definition~\InvLimIsolAtomsDef, and let
$\lambda_{\iota,n}:X_\iota\to X_{\iota,n}$ and
$\mu_{\iota,n}=(\lambda_{\iota,n})_*(\mu_\iota)$ the corresponding map and
measure.
Thus, we may regard $C(X_{\iota,n})$ as a unital C$^*$--subalgebra of
$C(X_\iota)$, where the state $\int\cdot\dif\mu_\iota$ restricts to
$\int\cdot\dif\mu_{\iota,n}$, and
$C(X)=\overline{\bigcup_{n=1}^\infty C(X_{\iota,n})}$.
If $y\in X_{\iota,n}$ is an atom of $\mu_{\iota,n}$ and if
$\lambda^{-1}_{\iota,n}(\{y\})$ contains no atoms of $\mu_\iota$ then since
$\lambda^{-1}_{\iota,n}(\{y\})$ is clopen in $X_\iota$ we may change
$(X_{\iota,n},\mu_{\iota,n})$ by substituting $\lambda^{-1}_{\iota,n}(\{y\})$
for $y$ and changing $\mu_{\iota,n}$ accordingly.
Then we must modify every $(X_{\iota,n+k},\mu_{\iota,n+k})$ too.
By doing so, we may assume without loss of generality that whenever
$n\in\Naturals$
and $y\in X_{\iota,n}$ is an atom of $\mu_{\iota,n}$ then
$\lambda_{\iota,n}^{-1}(\{y\})$ contains an atom of $\mu_\iota$.

Let $c=\max\{\mu_\iota(\{x\})\mid x\in X_\iota,\,1\le\iota\le N\}$.
Then $c<1$.
If $x_\iota\in X_\iota$ appears as one of the co--ordinates in an element of
$L\cup I_+$ then $1-\mu_\iota(\{x_\iota\})+(N-1)(1-c)\le1$, so
$\mu_\iota(\{x_\iota\})\ge(N-1)(1-c)$.
Therefore, $L\cup I_+$ is finite.
Moreover, if $\mu_\iota$ has infinitely many atoms then their masses form a
sequence tending to zero, so there is $\epsilon>0$ such that whenever
$x_\iota\in X_\iota$ $(1\le\iota\le N)$ and
$(x_\iota)_{\iota=1}^N\not\in L\cup I_+$ then
$\sum_{\iota=1}^N(1-\mu(\{x_\iota\}))>1+\epsilon$.

Let $x_\iota\in X_\iota$ be an atom of $\mu_\iota$.
If $x_\iota$ is an isolated point of $X_\iota$ then for $n$ large enough
$\lambda^{-1}_{\iota,n}(\{\lambda_{\iota,n}(x_\iota)\})=\{x_\iota\}$ and hence
the characteristic function of $\{x_\iota\}$, which we called $p_{x_\iota}$, is
in $C(X_{\iota,n})$.
We assume without loss of generality that this holds for every $n$ and for
every $\iota$.
If $x_\iota$ is not an isolated point of $X_\iota$ then
$\bigl(\lambda^{-1}_{\iota,n}
(\{\lambda_{\iota,n}(x_\iota)\})\bigr)_{n=1}^\infty$
is a neighborhood base for $x_\iota$ and, since $\mu_\iota$ is a regular
measure whose support is all of $X_\iota$, we have
$$ \mu_{\iota,n}(\{\lambda_{\iota,n}(x_\iota)\})
\ge\mu_{\iota,n+1}(\{\lambda_{\iota,n+1}(x_\iota)\})>\mu_\iota(\{x_\iota\}) $$
and
$\lim_{n\to\infty}\mu_{\iota,n}(\{\lambda_{\iota,n}(x_\iota)\})
=\mu_\iota(\{x_\iota\})$.
Thus, for $n$ large enough we have for every atom, $x_\iota\in X_\iota$ of
$\mu_\iota$, that
$$ \mu_{\iota,n}(\{\lambda_{\iota,n}(x_\iota)\})
<\mu_{\iota}(\{x_\iota\})+\epsilon/N. $$
We assume without loss of generality this holds for every $n$ and for every
$\iota$.

Then whenever $x_{\iota,n}\in X_{\iota,n}$ is an atom of $\mu_{\iota,n}$ and
$$ \sum_{\iota=1}^N(1-\mu_{\iota,n}(x_{\iota,n}))\le1, \tag{\lambdanatoms}  $$
there is a unique
$(x_\iota)_{\iota=1}^N\in L\cup I_+$ such that
$x_{\iota,n}=\lambda_{\iota,n}(x_\iota)$, $(1\le\iota\le N)$.
Moreover, if $(x_\iota)_{\iota=1}^N\in I_+$ then the atoms
$\lambda_{\iota,n}(x_\iota)$ and $x_\iota$ have the same mass, and if equality
holds in~(\lambdanatoms) then each $x_\iota$ is an isolated point of $X_\iota$.

Let $\Afr_n$ be the C$^*$--subalgebra of $\Afr$ generated by
$\bigcup_{\iota=1}^N C(X_{\iota,n})$.
Then
$$ \Afr_n\cong\freeprod{\iota=1}N(C(X_{\iota,n}),\tint\cdot\dif\mu_{\iota,n})
$$
and by Corollary~\FinManyCofXAtoms{} and the facts discussed above we have
$$ \Afr_n=\biggl(\smp{\Afr_{n,0}}{r_{n,0}}
\oplus\bigoplus_{x\in L_+}\smdp\Cpx{\alpha_{n,x}}{r_{n,x}}\biggr)
\oplus\bigoplus_{x\in I_+}\smdp\Cpx{\alpha_x}{r_x}, $$
where $I_+$, $\alpha_x$ and $r_x$ are as in the statement of the proposition
and where
$$ \align
L_+&=\left\{x=(x_\iota)_{\iota=1}^N\in L\bigg|
\dsize\sum_{\iota=1}^N(1-\mu_\iota(\{x_\iota\}))<1\right\} \\
\alpha_{n,x}&=1-\sum_{\iota=1}^N
(1-\mu_{\iota,n}(\{\lambda_{\iota,n}(x_\iota)\})) \\
r_{n,x}&=\bigwedge_{\iota=1}^N p_{x_\iota,n}
\endalign $$
where $p_{x_\iota,n}\in C(X_{\iota,n})$ is the characteristic function of
$\{\lambda_{\iota,n}(x_\iota)\}$.
Also, letting $L_0=L\backslash L_+$,
for each $x\in L_0$ there is a $*$--homomorphism
$\pi_{n,x}:\Afr_{n,0}\to\Cpx$ sending $p_{x_\iota,n}=p_{x_\iota}$ to $1$,
$(1\le\iota\le N)$.

We have $\Afr_n\subseteq\Afr_{n+1}$ and
$\Afr=\overline{\bigcup_{n=1}^\infty\Afr_n}$.
Let
$$ \Afr_{n,0}'=\smp{\Afr_{n,0}}{r_{n,0}}
\oplus\bigoplus_{x\in L_+}\smdp\Cpx{\alpha_{n,x}}{r_{n,x}}. $$
Then $\Afr_{n,0}'\subseteq\Afr_{n+1,0}'$ and~(\AfrIndLim) holds with
$\Afr_0=\overline{\bigcup_{n=1}^\infty\Afr_{n,0}'}$.
If $L=\emptyset$ then each $\Afr_{n,0}'=\Afr_{n,0}$ is simple with unique
tracial state and thus their inductive limit, $\Afr_0$, is simple with unique
tracial state.
Suppose $L$ is nonempty.
For each $x=(x_\iota)_{\iota=1}^N\in L$ let $\pi_{x,n}:\Afr_{n,0}'\to\Cpx$ be
the $*$--homomorphism sending $p_{x_\iota,n}\mapsto1$, $(1\le\iota\le N)$.
Then
$$ \Afr_{n,00}\eqdef\bigcap_{x\in L}\ker\pi_{x,n} $$
is simple with unique tracial state.
Clearly $\pi_{x,n+1}$ is an extension of $\pi_{x,n}$, so taking the inductive
limit we get, for each $x\in L$ a $*$--homomorphism $\pi_x:\Afr_0\to\Cpx$.
When restricted to $C(X_\iota)$, $\pi_x$ gives evaluation at
$x_\iota\in X_\iota$.
Then $\Afr_{00}$ is the inductive limit of the algebras $\Afr_{n,00}$, and thus
is simple with unique tracial state.

We now show that $\Afr_{00}$ is nonunital.
First consider the case when $L_+$ is nonempty.
Then for each $x=(x_\iota)_{\iota=1}^N\in L_+$, $\alpha_{n,x}$ decreases to the
limit $\sum_{\iota=1}^N(1-\mu_\iota(\{x_\iota\}))$ but is never equal to this
quantity.
Suppose for contradiction that $e\in\Afr_{00}$ is the identity of $\Afr_{00}$.
Let $n\in\Naturals$ and $a\in\Afr_{n,00}$.
Then there is $m>n$ such that $\alpha_{n,x}>\alpha_{m,x}$ and thus
$r_{n,x}-r_{m,x}\in\Afr_{m,00}$ is a nonzero projection.
Since $a\in\ker\pi_{n,x}$, since $\pi_{n,x}(r_{n,x})=1$ and since $r_{n,x}$ is
a minimal
projection of $\Afr_{n,00}$, we must have $ar_{n,x}=0$ and hence
$a(r_{n,x}-r_{m,x})=0$.
But $e(r_{n,x}-r_{m,x})=r_{n,x}-r_{m,x}$, so
$$ 1=\nm{r_{n,x}-r_{m,x}}=\nm{(e-a)(r_{n,x}-r_{m,x})}\le\nm{e-a}. $$
This contradicts that $e\in\overline{\bigcup_n\Afr_{n,00}}$.

If $L_+=\emptyset$ then $L_0$ is nonempty, so each $\Afr_{n,00}$ is nonunital.
But this implies that their inductive limit, $\Afr_{00}$, is nonunital.

The same technique shows that each ideal $\bigcap_{x\in F}\ker\pi_x$ of
$\Afr_0$ is nonunital.
\QED

Although the above proposition was stated only for free products of abelian
C$^*$--algebras, a similar result is easily proved for free products of
inductive limits of the algebras of the form $A_0\oplus\Cpx\cdots\oplus\Cpx$
that were considered, for example, in Proposition~\FinManyAzeroCdC.

\vskip3ex
\noindent{\bf\S\FPInfManyS. Free products of infinitely many algebras.}
\vskip3ex

In this section we consider the reduced free product of infinitely many
finite dimensional abelian C$^*$--algebras.
Although such free products of infinitely many  algebras can fail to be simple,
they never get a copy of $\Cpx$ as a direct summand, and hence their proper,
nontrivial ideals, if any, are always nonunital.
Moreover, the center of the free product algebra is always trivial,
even when its von Neumann algebra closure (i.e. the
strong--operator closure of the GNS representation) has projections that are
both minimal and central.

\proclaim{Theorem \FPInfMany}
For each $\iota\in\Naturals$ let $(A_\iota,\tau_\iota)$ be a finite dimensional
abelian algebra with faithful state as in~(\Aiota).
Let
$$ (\Afr,\tau)=\operatornamewithlimits*_{\iota=1}^\infty(A_\iota,\tau_\iota).
$$
Then $\Afr$ has stable rank~1.
Let
$$ L=\left\{\bigl(j(\iota)\bigr)_{\iota=1}^\infty\;\bigg|\;
\sum_{\iota=1}^\infty(1-\alpha_{\iota,j(\iota)})\le1\right\}. $$
If $L$ is empty then $\Afr$ is simple and $\tau$ is the unique
tracial state on $\Afr$.

Otherwise, if $L$ is nonempty, then for each
$j=(j(\iota))_{\iota=1}^\infty\in L$ there is a $*$--homomorphism
$\pi_j:\Afr\to\Cpx$ such that
$\pi_j(p_{\iota,j(\iota)})=1$ for every $\iota\in\Naturals$.
Then
$$ \Afr_{00}\eqdef\bigcap_{j\in L}\ker\pi_j $$
is simple and nonunital.
Moreover, letting
$$ \gamma_0=1-\sum_{j\in L}\left(
\sum_{\iota=1}^\infty\bigl(1-\alpha(\iota,j(\iota))\bigr)\right), $$
$\gamma_0^{-1}\tau\restrict_{\Afr_{00}}$ is the unique tracial state on
$\Afr_{00}$.
Finally, for every nonempty subset $F\subseteq L$, the ideal
$\bigcap_{j\in F}\ker\pi_j$ is nonunital.
\endproclaim
\demo{Proof}
The proof follows in a straightforward manner by using
Proposition~\FinManyFDAbel{} and taking inductive limits.
To show that $\Afr_{00}$ and each of the other proper, nontrivial ideals of
$\Afr$ are nonunital, an argument like that in the proof
Proposition~\InvLimIsolAtomsFP{} (in the case $L_+$ nonempty) is used.
\QED

Of course, similar results for free products of infinitely many abelian
algebras of the form considered in~\S\CofXS{} or infinitely many algebras of
the form $A_0\oplus\Cpx\oplus\cdots\oplus\Cpx$ considered in
Proposition~\FinManyAzeroCdC{} are easily obtained.

\vskip3ex
\noindent{\bf\S\Conjectures. Conjectures.}
\vskip3ex

This section contains a couple of related open problems which seem likely
to have solutions, though I don't yet see how to find them.

\proclaim{Conjecture~\ConjCofX}
Let $C(X)$ and $C(Y)$ be unital, abelian C$^*$--algebras having faithful states
given by probability measures $\mu_X$ on $X$ and $\mu_Y$ on $Y$.
Let
$$ (\Afr,\tau)=(C(X),\tint\cdot\dif\mu_X)*(C(Y),\tint\cdot\dif\mu_Y). $$
Then necessary and sufficient conditions for $\Afr$ to be simple is that for
every $x\in X$ and $y\in Y$, $\mu_X(\{x\})+\mu_Y(\{y\})<1$.
\endproclaim

\proclaim{Proposition \CofXNec}
The conditions named in Conjecture~\ConjCofX{} are necesary for
simplicity of $\Afr$.
\endproclaim
\demo{Proof}
Suppose the conditions of~\ConjCofX{} are not satisfied.
Let $C\subseteq X$ be the set of atoms of $\mu_X$.
If $C$ is finite then let $D=C$.
If $C$ is infinite then let $\overline C$ denote the closure of $C$ in $X$ and
let $D$ be a totally disconnected, separable, compact, Hausdorff
space equipped with a continuous, surjective map $\lambda:D\to\overline C$.
(That such a space exists is a well--known result.)
For each $c\in C$, choose $f(c)\in\lambda^{-1}(D)$.
Let $\nu_{D}$ be the measure on $D$ given by
$\nu_{D}(\{f(c)\})=\mu_X(\{c\})$ and $\nu_{D}(D\backslash f(C))=0$.
Replace $D$ by the support of $\nu_D$.
Let $\nu_{X'}$ be the measure on $X$ that when restricted to $X\backslash C$
gives the same measure as $\mu_X$ and such that $\nu_{X'}(C)=0$ and let $X'$ be
the support of $\nu_{X'}$.
Let $X_a$ be the compact Hausdorff space that is the disjoint union of $X'$ and
$D$ and let $\mu_{X_a}$ be the measure on $X_a$ obtained from $\nu_{X'}$ and
$\nu_D$.
Let $\kappa_X:X_a\to X$ be the surjective, continous map composed of
the inclusion $X'\hookrightarrow X$ and $\lambda:D\to X$.
Then $\mu_{X_a}$ is the push--forward measure,
$\mu_{X_a}=(\kappa_X)_*(\mu_X)$.
Therefore $\kappa_X$ induces an injective, unital $*$--homomorphism,
$\pi_X:C(X)\to C(X_a)$, preserving the states defined by the measures
$\mu_{X_a}$ and $\mu_X$.
Do the same for $Y$, getting $\pi_Y:C(Y)\to C(Y_A)$.

Now Proposition~\InvLimIsolAtomsFP{} applies to the free product
$$ (\Afr',\tau')=(C(X_a),\tint\cdot\dif\mu_{X_a})
*(C(Y_a),\tint\cdot\dif\mu_{Y_a}). $$
If condition of Conjecture~\ConjCofX{} is not satisfied then there are
$x\in X_a$ and $y\in Y_a$ such that
$\mu_{X_a}(\{x\})+\mu_{Y_a}(\{y\})\ge1$, which by~\InvLimIsolAtomsFP{} implies
the existence of a $*$--homomorphism $\pi:\Afr'\to\Cpx$ that when restricted to
$C(X_a)$ is evaluation at $x$.
Since $\tau'$ is faithful and $\kappa_X$ and $\kappa_Y$ are trace--preserving
inclusions,
they induce an inclusion $\Afr\hookrightarrow\Afr'$.
Then $\pi\restrict_{\Afr}$ is a nonzero $*$--homomorphism, so $\Afr$ is not
simple.
\QED

We now look at free products of finite dimensional algebras, and
we use the notation of~\cite{\DykemaZZFPFDNT} for
a faithful state on a finite dimensional algebra.
Thus, if $D=\bigoplus_{j=1}^K M_{n_j}(\Cpx)$, we write
$$ (D,\phi)=\bigoplus_{j=1}^K \smd{M_{n_j}(\Cpx)}
{\alpha_{j,1},\ldots,\alpha_{j,n_j}} $$
to mean that the restriction of $\phi$ to the $j$th summand of $D$ is given by
$\Tr(\cdot H)$ where $\Tr$ is the trace on $M_{n_j}(\Cpx)$ of norm $n_j$ and
where $H$ is the
diagonal matrix with $\alpha_{j,1},\ldots,\alpha_{j,n_j}$ down the
diagonal.
\proclaim{Conjecture~\ConjFD}
Let
$$ \aligned
(A_1,\phi_1)&=\bigoplus_{j=1}^{K_1}\smd{M_{n_j}(\Cpx)}
{\alpha_{j,1},\cdots,\alpha_{j,n_j}} \\
(A_2,\phi_2)&=\bigoplus_{j=1}^{K_2}\smd{M_{m_j}(\Cpx)}
{\beta_{j,1},\cdots,\beta_{j,m_j}}
\endaligned  $$
be finite dimensional C$^*$--algebras with faithful states and let
$$ (\Afr,\phi)=(A_1,\phi_1)*(A_2,\phi_2)   \tag{\ConjFDFP} $$
be the reduced free product C$^*$--algebra.
Then necessary and sufficient conditions for $\Afr$ to be simple are that
if $n_j=1$ for some $j$ then for every $1\le k\le K_2$
$$ \frac1{1-\alpha_{j,1}}<\sum_{i=1}^{m_k}\frac1{\beta_{k,i}} $$
and if $m_j=1$ for some $j$ then for every $1\le k\le K_1$
$$ \frac1{1-\beta_{j,1}}<\sum_{i=1}^{n_k}\frac1{\alpha_{k,i}}. $$
\endproclaim
If the above conjecture is true, then in particular $\Afr$ is simple whenever
each $n_j>1$ and each $m_j>1$.
This conjecture is inspired by the results of~\cite{\DykemaZZFreeDim}
and~\cite{\DykemaZZFPFDNT}, which
show that necessary and sufficient conditions for the von Neumann algebra free
product analogous to~(\ConjFDFP) to be a
factor are that
if $n_j=1$ for some $j$ then for every $1\le k\le K_2$
$$ \frac1{1-\alpha_{j,1}}\le\sum_{i=1}^{m_k}\frac1{\beta_{k,i}}
\tag{\ConjFDnj} $$
and if $m_j=1$ for some $j$ then for every $1\le k\le K_1$
$$ \frac1{1-\beta_{j,1}}\le\sum_{i=1}^{n_k}\frac1{\alpha_{k,i}}.
\tag{\ConjFDmj} $$
Moreover, using this von Neumann algebra result and an argument similar to the
one used in Proposition~\CofXNec{} we see that~(\ConjFDnj) and~(\ConjFDmj) are
necessary conditions for the simplicity of $\Afr$.

For some results about certain reduced free product C$^*$--algebras with
respect to non--tracial states, see~\cite{\DykemaRordamZZPI}.

\enddocument